\documentclass[11pt]{article}

\usepackage{preprint-layout-light}
\usepackage{mathtools}
\usepackage{bm}
\usepackage{booktabs}
\usepackage{placeins}
\usepackage{pgfplots}

\pgfplotsset{compat=1.18,log basis x=10}
\graphicspath{{figures/}}

\newcommand{\R}{\mathbb{R}}
\newcommand{\Pp}{\mathbb{P}}
\newcommand{\Th}{\mathcal{T}_h}
\newcommand{\Tzero}{\mathcal{T}_{0,h}}

\newcommand{\Ph}{\mathcal{P}_h}

\newcommand{\energy}[1]{\mathopen{\vert\!\vert\!\vert}#1\mathclose{\vert\!\vert\!\vert}}
\newcommand{\normal}{n}

\newcommand{\addressumushort}{Mathematics and Mathematical Statistics, Ume{\aa}~University, Sweden}
\newcommand{\addressjushort}{Mechanical Engineering, J\"onk\"oping University, Sweden}
\newcommand{\addressuclshort}{Mathematics, University College London, UK}

\title{One-Sided $p$-Explicit Stability for High-Order\\ Cut Finite Element Nitsche Methods}
\author{Erik Burman \and Peter Hansbo \and Mats G. Larson \and Karl Larsson}
\date{}

\begin{document}

\maketitle

\begin{abstract}
We give a one-sided stability argument for high-order cut finite elements with symmetric Nitsche boundary conditions. On each boundary patch, an $H^1$-projection splits the discrete function into a single polynomial and a stabilized residual. Inward segments control the polynomial inside the physical domain; complete-element estimates control only the residual. For smooth exterior boundaries and under a geometry-dependent mesh restriction, this gives an explicit normal-flux estimate with factor $p^2/h$ and a sufficient Nitsche penalty with the same scaling. The constants are independent of the cut position and polynomial degree. A neighboring-interval estimate quantifies the exponential cost of polynomial continuation and motivates avoiding it in the stability proof. We obtain $hp$-error estimates, with half a power of degree lost in energy and a full power in $L^2$ relative to best approximation, and exponential degree convergence for analytic extensions on sufficiently fine meshes. Translated-disk experiments show bounded normalized flux constants and penalty thresholds over the tested degrees and cuts, and resolved convergence curves follow the predicted orders. Physical and extension-dependent stabilization errors are reported separately.
\end{abstract}

\section{Introduction}
\label{sec:introduction}

Unfitted finite element methods permit the physical boundary to intersect a
background mesh. Cut finite element methods (CutFEM) combine this flexibility
with weak formulations and stabilization near the cut boundary
\cite{Burman2010,BurmanHansbo2012,BurmanEtAl2015,CutFEM2025}.
Symmetric Nitsche stability requires control of the discrete normal flux.
On fitted elements, a polynomial inverse trace inequality gives the penalty
scale \(p^2/h\). On cut elements, arbitrarily small physical intersections
prevent a direct use of that argument.

For interface-fitted subdomain meshes that need not match across the interface,
Chernov and Hansbo~\cite{ChernovHansbo2011} developed an \(hp\)-Nitsche method
with degree-dependent penalties and error estimates. High-order unfitted
methods have been developed along several complementary directions. These include \(hp\)-explicit interface penalty methods
\cite{WuXiao2019}, high-order geometry approximation by isoparametric mappings
\cite{Lehrenfeld2016,LehrenfeldReusken2018}, higher-order ghost penalties
\cite{StickoKreiss2019,Wichrowski2026}, unfitted spectral element methods
\cite{GonzalezGuoYang2025}, and cell-merging high-order interface methods with
explicit conditioning analysis \cite{ChenLiu2023,JohanssonLarson2013}. These works establish that
high-order and spectral accuracy are compatible with unfitted geometry. The question
addressed here is more specific: how should one prove symmetric Nitsche stability
for an exterior boundary with constants that remain explicit and polynomial in \(p\)?

Standard CutFEM analysis first extends control from the physical domain to
complete active elements and then applies a full-element inverse trace
inequality. This is effective for fixed degree and supplies active-mesh
control for condition-number estimates. Aggregation and extended spaces
modify the trial space \cite{BadiaVerdugoMartin2018,BurmanHansboLarson2022};
ghost penalties enforce control weakly \cite{Burman2010,BadiaNeivaVerdugo2022Linking,
BurmanHansboLarson2026}. Macroelements combining large and small physical
intersections also permit locally conservative discontinuous CutFEM
\cite{LarsonZahedi2023}.

The degree dependence of full-element transfer matters both for particular
extension operators and for proofs that use them. Extrapolation-based high-order
aggregation has exponentially growing stability constants, whereas an
interpolation-based discrete extension restores polynomial high-order robustness
and fitted-like conditioning \cite{BadiaNeivaVerdugo2022HighOrder}.
Section~\ref{sec:continuation-barrier} records the sharp neighboring-interval
\(L^2\) continuation order as motivation. An exponential norm of that
continuation operator does not imply that the underlying Nitsche method needs
an exponential penalty: the intermediate full-element estimate may be stronger
than coercivity requires.

Our argument instead projects onto one polynomial per boundary patch.
Inward normal segments control that polynomial inside the physical domain;
a gradient-projection stabilization controls the residual on complete elements.
Here \emph{one-sided} refers to the polynomial's physical-side control; the
Nitsche form remains symmetric.

Physical-domain coercivity does not imply a cut- or \(p\)-independent matrix
condition number: small cuts can still create nearly dependent basis functions
\cite{dePrenterEtAl2017,dePrenterEtAl2023}. Aggregation and algebraic or subspace
preconditioners address this separate issue
\cite{BadiaVerdugoMartin2018,BadiaVerdugo2018,GrossReusken2023}.
Our numerical patch Schwarz comparison reports its remaining degree dependence.

\paragraph{Main Results.}
The main results concern stability and its approximation consequences.

\begin{itemize}
\item On the constructed boundary patches, we prove the one-sided inverse estimate
\begin{equation}
h\|\partial_n v_h\|_{\partial\Omega}^2
\lesssim
p^2\bigl(\|\nabla v_h\|_\Omega^2+\|v_h\|_{s_h}^2\bigr)
\qquad \forall v_h\in V_h \label{eq:introduction-one-sided}
\end{equation}
with the explicit constant in \eqref{eq:assembled-inverse-constant}, independent
of the cut configuration, \(h\), and \(p\). For \(h\le h_0\), with \(h_0\)
determined by the boundary geometry and mesh regularity, it yields symmetric
Nitsche coercivity with penalty coefficient \(\beta_0p^2/h\) above an explicit,
cut- and degree-independent sufficient value of \(\beta_0\).
\item Under the stated regularity assumptions, the same stabilization gives a
Strang bound, \(hp\)-energy and \(L^2\)-error estimates, and exponential degree
convergence for analytically extendable solutions on sufficiently fine meshes.
Relative to optimal approximation, the energy and \(L^2\) bounds lose half a
power and a full power of \(p\), respectively.
\end{itemize}

The contribution is the stability argument for this projection-stabilized
symmetric formulation. Translated-disk experiments test its flux and convergence
scaling. Their meshes do not meet the sufficient restriction in the explicit
disk calculation; measured thresholds are distinguished from analytic bounds.

\paragraph{Position Relative to the Closest Work.}
The one-sided idea has an important interface analogue. Wu and Xiao
\cite{WuXiao2019} decompose an arithmetic interface flux into the flux from either
side plus a controlled flux jump, so that an inverse estimate is needed on only one
of the two physical subcells. An exterior boundary has no complementary physical
subcell. The patch projection and inward-segment argument introduced here replace
that missing side while preserving the same \(p^2/h\) scaling.

Projection-type ghost penalties and their relation to derivative-jump penalties and
discrete extensions are established in
\cite{BurmanHansbo2014Stokes,BadiaNeivaVerdugo2022Linking,
BurmanHansboLarson2026}. Our contribution is the \(p\)-explicit one-sided flux
estimate: complete-element inverse estimates are applied only to the projection
residual, while the projected polynomial is controlled inside the physical domain.
This preserves the conforming active-mesh space and symmetric Nitsche formulation
with a fitted-like penalty scale.

\paragraph{Organization.}
Section \ref{sec:setting} defines the method and constructs the boundary patches.
Section \ref{sec:continuation-barrier} identifies the exponential continuation
barrier. Section \ref{sec:one-sided-stability} proves the one-sided inverse estimate
and Nitsche coercivity. Section \ref{sec:error-analysis} gives the \(hp\)-error
analysis. Section \ref{sec:numerical-results} presents the numerical studies, and
Section \ref{sec:conclusions} summarizes the conclusions.

\section{Problem Setting and Discrete Method}
\label{sec:setting}

\paragraph{Standing Assumptions.}
The following assumptions are used throughout the paper.
\begin{enumerate}[label=(A\arabic*),leftmargin=*,itemsep=2pt,topsep=3pt]
\item The physical domain \(\Omega\subset\R^d\), \(d\in\{2,3\}\), is bounded and
  has a \(C^2\) boundary. The volume and boundary integrals in the method are taken
  over the exact sets \(\Omega\) and \(\partial\Omega\); geometry-approximation and
  quadrature errors are not included in the analysis.
\item The background meshes are shape-regular, quasi-uniform affine simplicial
  meshes, and the polynomial degree \(p\ge1\) is uniform on each mesh.
\item The mesh size satisfies \(0<h\le h_0\), where the geometry-dependent threshold
  \(h_0\) is fixed in Proposition~\ref{prop:boundary-patch-covering} and may be
  reduced, independently of \(p\) and the cut position, in subsequent local geometric
  arguments. All patch-construction parameters are fixed independently of \(h\),
  \(p\), and the cut position.
\item Unless displayed explicitly, constants may depend only on the dimension, the
  \(C^2\)-geometry of \(\partial\Omega\), the mesh shape-regularity and
  quasi-uniformity constants, and the fixed patch parameters. They are independent of
  \(h\), \(p\), and how \(\partial\Omega\) cuts the mesh.
\end{enumerate}

\subsection{Model Problem}
\label{subsec:model-problem}

Let \(\normal\) denote the exterior unit normal to \(\partial\Omega\). We consider
the Dirichlet problem
\begin{alignat}{2}
-\Delta u &= f &\qquad&\text{in }\Omega \label{eq:model-bulk}\\
u &= g &&\text{on }\partial\Omega \label{eq:model-boundary}
\end{alignat}
where \(f\in L^2(\Omega)\) and \(g\in H^{1/2}(\partial\Omega)\).
For consistency and error estimates involving the classical normal trace,
we additionally assume \(u\in H^{3/2+\varepsilon}(\Omega)\) for some
\(\varepsilon>0\), so that \(\partial_nu\in L^2(\partial\Omega)\).
The approximation theorems below impose their stated stronger regularity.

For sufficiently small \(\delta_0>0\), the inward normal-coordinate map is defined by
\begin{equation}
\Phi:\partial\Omega\times(0,\delta_0)\rightarrow\Omega,
\qquad
\Phi(y,t)=y-t\normal(y) \label{eq:normal-coordinate-map}
\end{equation}
and is a diffeomorphism onto the interior tubular neighborhood \(U_{\delta_0}\).
Its Jacobian satisfies
\begin{equation}
J(y,t)=\prod_{i=1}^{d-1}\bigl(1-t\kappa_i(y)\bigr),
\qquad
0<c_J\le J(y,t)\le C_J \label{eq:normal-coordinate-jacobian}
\end{equation}
after reducing \(\delta_0\) if necessary.

\subsection{Active Mesh and Finite Element Space}
\label{subsec:active-mesh}

Let \(\Omega_0\) be a polygonal or polyhedral domain containing
\(\overline\Omega\). Let \(\Tzero\) be a shape-regular, quasi-uniform family of
affine simplicial meshes of \(\Omega_0\), with mesh parameter \(h>0\).
We use fixed constants \(C_T,\varrho>0\) such that
\(\operatorname{diam}(T)\le C_T h\) and every element contains a ball
of radius \(\varrho h\). Thus \(h\) need not equal the maximal element
diameter; in the right-triangular numerical mesh it is the square spacing
and \(C_T=\sqrt2\).
The active mesh and its domain are
\begin{equation}
\Th=\{T\in\Tzero:T\cap\Omega\ne\varnothing\},
\qquad
\Omega_h=\bigcup_{T\in\Th}T \label{eq:active-mesh}
\end{equation}
and the cut elements are
\begin{equation}
\mathcal T_{\partial,h}
=
\{T\in\Th:T\cap\partial\Omega\ne\varnothing\} \label{eq:cut-elements}
\end{equation}

For a fixed polynomial degree \(p\ge1\), define
\begin{equation}
V_{0,h}
=
\{v\in C^0(\Omega_0):v|_T\in\Pp_p(T)\ \forall T\in\Tzero\},
\qquad
V_h=V_{0,h}|_{\Omega_h} \label{eq:finite-element-space}
\end{equation}
where \(\Pp_p(T)\) denotes polynomials of total degree at most \(p\).

\subsection{Constructed Boundary Patches}
\label{subsec:boundary-patches}

The following covering result supplies the patches used by the method. It is a
consequence of the tubular-neighborhood theorem and the local counting property of a
quasi-uniform mesh.

\begin{figure}[htbp]
\centering
\resizebox{0.65\linewidth}{!}{\begingroup
  \fontsize{16}{16.8}\selectfont
  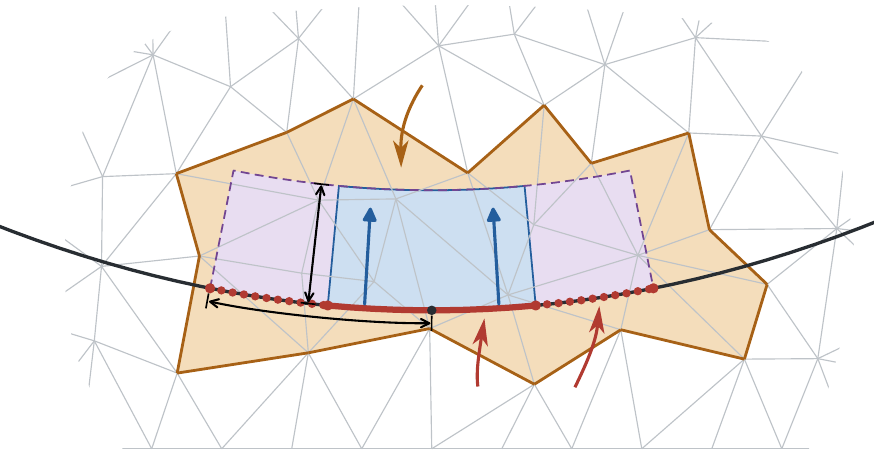\endgroup
}
\caption{Boundary-patch construction.
Choose \(y_i\) from a maximal \(h\)-separated set on \(\partial\Omega\) and take
its geodesic Voronoi cell \(\Gamma_{\mathcal P_i}\) (solid red).
Enlarge the cell to \(\omega_i=B_\Gamma(y_i,R_\Gamma h)\) (solid and dotted red),
with \(R_\Gamma>1\) and radius measured along \(\partial\Omega\).
Sweep the cell and its enlargement inward along normals through depth
\(c_\delta h\) to form \(U_{\mathcal P_i}\) (blue) and
\(W_{\mathcal P_i}\) (purple), respectively.
Collect all complete active elements whose interiors meet \(W_{\mathcal P_i}\)
to form the face-connected patch \(\mathcal P_i\) (brown).}
\label{fig:boundary-patch-construction}
\end{figure}

Figure~\ref{fig:boundary-patch-construction} illustrates how a boundary cell and
its enlarged inward tube determine a complete-element patch.

\begin{prop}[Boundary-Patch Covering]
\label{prop:boundary-patch-covering}
There are constants \(h_0,c_\delta,c_{\mathcal P},C_{\mathcal P}>0\) and
\(N_{\mathcal P}\in\mathbb N\), independent of \(h\) and of how
\(\partial\Omega\) cuts the mesh, with the following property. For every
\(h\in(0,h_0]\), one can construct a partition
\(\{\Gamma_{\mathcal P}\}_{\mathcal P\in\Ph}\) of \(\partial\Omega\) and
face-connected patches \(\mathcal P\subset\Omega_h\), each formed as the union of
complete active background elements, such that
\begin{align}
c_{\mathcal P}h
&\le
\operatorname{diam}(\mathcal P)
\le
C_{\mathcal P}h \label{eq:patch-diameter}\\
U_{\mathcal P}
&:=
\{\Phi(y,t):y\in\Gamma_{\mathcal P},\ 0<t<c_\delta h\}
\subset
\mathcal P\cap\Omega \label{eq:patch-inward-tube}\\
\sum_{\mathcal P\in\Ph}\mathbf{1}_{\mathcal P}(x)
&\le
N_{\mathcal P}
\qquad \forall x\in\Omega_h \label{eq:patch-overlap}
\end{align}
The boundary pieces are pairwise disjoint up to sets of surface measure zero and
satisfy
\begin{equation}
\partial\Omega
=
\bigcup_{\mathcal P\in\Ph}\overline{\Gamma}_{\mathcal P} \label{eq:boundary-partition}
\end{equation}
Each patch contains a uniformly bounded number of elements.
\end{prop}

\begin{proof}
Let \(d_\Gamma\) denote geodesic distance on \(\partial\Omega\). Compactness and
\(C^2\)-regularity provide a radius \(r_\Gamma>0\) and a constant
\(C_d\ge1\) on which geodesic and Euclidean distances are \(C_d\)-equivalent,
geodesic balls are coordinate neighborhoods, and
their surface measures satisfy
\begin{equation}
c_\Gamma r^{d-1}
\le
|B_\Gamma(y,r)|
\le
C_\Gamma r^{d-1},
\qquad 0<r\le r_\Gamma 
\label{eq:surface-ball-regularity}
\end{equation}
We take \(r_\Gamma\) smaller than a fixed fraction of the reciprocal maximum
curvature. Fix \(c_\delta>0\) and a surface-enlargement factor \(R_\Gamma>1\), set
\(C_0=R_\Gamma+c_\delta+C_T\), and choose
\begin{equation}
h_0
\le
\min\left\{
\frac{\delta_0}{2c_\delta},
\frac{r_\Gamma}{4C_d(C_0+R_\Gamma+1)}
\right\}
\label{eq:patch-h0-choice}
\end{equation}
Thus every normal segment used below lies inside the tubular neighborhood, and every
surface neighborhood used in the construction lies in a single uniformly regular
chart.

For \(0<h\le h_0\), choose a maximal \(h\)-separated set
\(\{y_i\}_{i=1}^{M_h}\subset\partial\Omega\), with associated patches
\(\mathcal P_i\) to be constructed below. Then
\begin{equation}
d_\Gamma(y_i,y_j)\ge h\quad(i\ne j),
\qquad
\partial\Omega\subset\bigcup_{i=1}^{M_h}B_\Gamma(y_i,h)
\label{eq:maximal-boundary-net}
\end{equation}
Let \(\Gamma_{\mathcal P_i}\) be the geodesic Voronoi cell of \(y_i\), with ties assigned by a
fixed ordering. Then the cells form a measurable partition up to their boundaries and
\begin{equation}
B_\Gamma(y_i,h/2)\subset\overline{\Gamma_{\mathcal P_i}}
\subset\overline{B_\Gamma(y_i,h)}
\label{eq:voronoi-cell-bounds}
\end{equation}
Indeed, the first inclusion follows from the separation in
\eqref{eq:maximal-boundary-net}, and the second from maximality. Hence each cell has
diameter comparable with \(h\). Moreover, the enlarged neighborhoods
\(B_\Gamma(y_i,R_\Gamma h)\) have uniformly bounded overlap: the balls
\(B_\Gamma(y_i,h/3)\) are disjoint, and the surface-measure bounds
\eqref{eq:surface-ball-regularity} applied in a ball of radius
\((R_\Gamma+1)h\) give a bound depending only on \(d\), \(R_\Gamma\), and the
geometry constants \(c_\Gamma,C_\Gamma\).

Set \(\omega_i=B_\Gamma(y_i,R_\Gamma h)\). The inward tubes generated by
\(\Gamma_{\mathcal P_i}\) and \(\omega_i\), respectively, are
\begin{align}
U_{\mathcal P_i}&:=\{\Phi(y,t):y\in\Gamma_{\mathcal P_i},\ 0<t<c_\delta h\}
\label{eq:cell-inward-tube}\\
W_{\mathcal P_i}&:=\{\Phi(y,t):y\in\omega_i,\ 0<t<c_\delta h\}
\label{eq:enlarged-inward-tube}
\end{align}
Define \(\mathcal P_i\) as the union of all active elements whose interiors meet
\(W_{\mathcal P_i}\), as illustrated in Figure~\ref{fig:boundary-patch-construction}.
Because \(W_{\mathcal P_i}\) is connected, a path in \(W_{\mathcal P_i}\), perturbed off the mesh
skeleton of codimension two, shows that these elements are face-connected. The cells
satisfy \(\overline{\Gamma_{\mathcal P_i}}\subset\omega_i\) by
\eqref{eq:voronoi-cell-bounds} and \(R_\Gamma>1\), and therefore
\begin{equation}
U_{\mathcal P_i}\subset W_{\mathcal P_i}\subset\mathcal P_i\cap\Omega
\label{eq:cell-tube-in-patch}
\end{equation}

Using distance equivalence and quasi-uniformity, we obtain
\begin{equation}
\mathcal P_i\subset B(y_i,C_0h)
\label{eq:patch-ball-containment}
\end{equation}
where \(C_0=R_\Gamma+c_\delta+C_T\) is fixed. Every element has volume at least
\(c_Th^d\), while the element interiors are disjoint and all elements in
\(\mathcal P_i\) lie in \(B(y_i,(C_0+1)h)\). Consequently,
\begin{equation}
\#\{T:T\subset\mathcal P_i\}
\le
\frac{|B(y_i,(C_0+1)h)|}{c_Th^d}
\le N_T 
\label{eq:patch-element-count}
\end{equation}
This proves the uniform element count and the upper diameter bound in
\eqref{eq:patch-diameter}; the inward segment based at \(y_i\) gives the lower
diameter bound \(\operatorname{diam}(\mathcal P_i)\ge c_\delta h\).

Finally, if \(x\in\mathcal P_i\), then \eqref{eq:patch-ball-containment} places
\(y_i\) within \(C_0h\) of \(x\). The centers associated with all patches containing
the same point are therefore contained in a surface ball of radius \(Ch\). Their
surface balls of radius \(h/3\) are disjoint, so a second application of
\eqref{eq:surface-ball-regularity} bounds their number by a uniform
\(N_{\mathcal P}\). This is \eqref{eq:patch-overlap}. The constants
\(N_T,N_{\mathcal P},c_{\mathcal P}\), and \(C_{\mathcal P}\) depend only on the
dimension, the mesh shape-regularity and quasi-uniformity constants, the fixed choices
\(c_\delta,R_\Gamma\), and the \(C^2\)-geometry through
\(\delta_0,r_\Gamma,c_\Gamma,C_\Gamma\); they do not depend on the cut position.
\end{proof}

\begin{rem}[Checking Numerical Patches]
\label{rem:patch-certification}
The construction includes the inward tube by definition. Since each
face-connected \(\mathcal P_i\) contains at most \(N_T\) elements, the
active-mesh patch within \(N_T-1\) face steps of any seed in \(\mathcal P_i\)
contains the whole patch and its tube. This is a uniform sufficient depth,
not a minimal layer count. Smaller patches require a verified assignment
of boundary pieces to patches containing their tubes. For disks,
circle--edge and radial vertex events give a finite coverage check;
Appendix~\ref{subsec:geometry-audit} verifies assignments without enlargement.
Increasing \(c_\delta\) improves the flux constant but may require larger
patches, increasing the region approximated by one polynomial.
\end{rem}

\subsection{Patch Projection and Stabilization}
\label{subsec:patch-stabilization}

Let \(\Pp_p(\mathcal P)\) denote restrictions to \(\mathcal P\) of single global
polynomials in \(\Pp_p(\R^d)\). Define the \(H^1\)-projection
\(E_{\mathcal P}:H^1(\mathcal P)\rightarrow\Pp_p(\mathcal P)\) by
\begin{align}
(\nabla E_{\mathcal P}v,\nabla q)_{\mathcal P}
&=
(\nabla v,\nabla q)_{\mathcal P}
\qquad \forall q\in\Pp_p(\mathcal P) \label{eq:patch-projection-gradient}\\
(E_{\mathcal P}v,1)_{\mathcal P}
&=
(v,1)_{\mathcal P} \label{eq:patch-projection-mean}
\end{align}
The gradient form is positive definite on \(\Pp_p(\mathcal P)/\R\),
because the patch has nonempty interior, and the mean condition fixes the
constant. Thus the projection is well defined, reproduces every polynomial in
\(\Pp_p(\mathcal P)\), and is stable in the \(H^1\)-seminorm:
\begin{equation}
\|\nabla E_{\mathcal P}v\|_{\mathcal P}
\le
\|\nabla v\|_{\mathcal P} \label{eq:patch-projection-stability}
\end{equation}

The patch stabilization is
\begin{equation}
s_h(v,w)
=
\sum_{\mathcal P\in\Ph}
\bigl(\nabla(v-E_{\mathcal P}v),
\nabla(w-E_{\mathcal P}w)\bigr)_{\mathcal P} \label{eq:patch-stabilization}
\end{equation}
The associated stabilization seminorm is defined by
\begin{equation}
\|v\|_{s_h}^2
:=
\sum_{\mathcal P\in\Ph}
\|\nabla(v-E_{\mathcal P}v)\|_{\mathcal P}^2 
\label{eq:stabilization-seminorm}
\end{equation}

\begin{rem}[Stabilization Seminorm]
\label{rem:stabilization-seminorm}
The double bars in \(\|v\|_{s_h}\) denote a seminorm, not a norm. Indeed, the
quantity vanishes for every global polynomial in \(\Pp_p(\R^d)\) restricted to
\(\Omega_h\), because the patch projections reproduce such polynomials. The notation
\(\|\cdot\|_{s_h}\) is used throughout whenever the stabilization is evaluated on
the same function in both slots.
\end{rem}

\begin{rem}[Mean Preservation]
\label{rem:no-l2-stabilization}
The mean condition \eqref{eq:patch-projection-mean} and the uniform Poincar\'e
inequality on the face-connected complete-element patches give
\begin{equation}
h^{-2}\sum_{\mathcal P\in\Ph}\|v-E_{\mathcal P}v\|_{\mathcal P}^2
\lesssim \|v\|_{s_h}^2
\qquad \forall v\in H^1(\Omega_h)
\label{eq:l2-residual-control}
\end{equation}
The constant is independent of \(h\), \(p\), and the cut configuration:
element Poincar\'e inequalities and comparisons of neighboring element means
along a face-connected tree involve only a uniformly bounded number of
shape-regular, quasi-uniform elements. Thus a separate \(L^2\)-residual
penalty is unnecessary.
\end{rem}

In the analysis, the stabilization \(s_h\) enters with coefficient one.
Section~\ref{sec:numerical-results} discusses rescaling by a fixed positive
coefficient \(\gamma_1\).

\subsection{Nitsche Formulation}
\label{subsec:nitsche-formulation}

For piecewise smooth functions, \(\partial_n v\) on \(\partial\Omega\)
is understood almost everywhere as the trace from the physical side. This
convention also applies when the boundary coincides with a mesh face on a set
of positive surface measure. For sufficiently regular \(v,w\), define
\begin{align}
a_h(v,w)
&=
(\nabla v,\nabla w)_\Omega
-
(\partial_n v,w)_{\partial\Omega}
-
(v,\partial_n w)_{\partial\Omega}
+
\beta_0\frac{p^2}{h}(v,w)_{\partial\Omega} \label{eq:nitsche-form}\\
A_h(v,w)
&=
a_h(v,w)+s_h(v,w) \label{eq:stabilized-form}
\end{align}
and
\begin{equation}
l_h(w)
=
(f,w)_\Omega
-
(g,\partial_n w)_{\partial\Omega}
+
\beta_0\frac{p^2}{h}(g,w)_{\partial\Omega} \label{eq:linear-functional}
\end{equation}
The discrete method is: find \(u_h\in V_h\) such that
\begin{equation}
A_h(u_h,v_h)=l_h(v_h)
\qquad \forall v_h\in V_h \label{eq:discrete-method}
\end{equation}

The stabilization in \eqref{eq:patch-stabilization} is weakly consistent rather than
strongly consistent. Its contribution is retained explicitly in the Strang estimate
and bounded using polynomial reproduction and local approximation.

\section{The Full-Element Continuation Barrier}
\label{sec:continuation-barrier}

The standard ghost-penalty analysis obtains control on \(\Omega_h\) before applying
inverse inequalities. At high polynomial degree, the transfer between neighboring
complete elements is already exponentially expensive. The following elementary result
isolates this mechanism. The exponential base is the classical Chebyshev
growth factor outside the reference interval \cite[Eq.~18.5.1]{DLMF};
the proposition records its sharp \(L^2\) continuation order, including
the algebraic prefactor.

\begin{prop}[Sharp Order of Neighboring-Interval Continuation]
\label{prop:interval-continuation}
Let \(I_-=(-1,0)\), \(I_+=(0,1)\), and
\begin{equation}
C_p=\sup_{q\in\Pp_p(\R)\setminus\{0\}}
\frac{\|q\|_{I_+}}{\|q\|_{I_-}}
\label{eq:continuation-constant}
\end{equation}
Set \(\rho=3+2\sqrt2\) and \(a=\operatorname{arccosh}3=\log\rho\). For \(p\ge1\),
\begin{align}
C_p&\ge B_p:=\left(\frac{I_p}{D_p}\right)^{1/2},
\qquad D_p=1-\frac{1}{4p^2-1}
\label{eq:continuation-lower-bound}\\
I_p&=1+\frac14\left\{
\frac{\cosh((2p+1)a)-1}{2p+1}
-\frac{\cosh((2p-1)a)-1}{2p-1}\right\}
\label{eq:chebyshev-continuation-integral}
\end{align}
Moreover,
\begin{equation}
B_p\sim 2^{-3/4}p^{-1/2}\rho^p,
\qquad C_p\asymp p^{-1/2}\rho^p
\label{eq:continuation-sharp-order}
\end{equation}
The constants implicit in \(\asymp\) are independent of \(p\).
\end{prop}

\begin{proof}
Choose \(q(x)=T_p(2x+1)\), where \(T_p\) is the Chebyshev polynomial.
The change of variable \(z=2x+1\) gives
\begin{equation}
\frac{\|q\|_{I_+}^2}{\|q\|_{I_-}^2}
=\frac{\int_1^3T_p(z)^2\,dz}{\int_{-1}^1T_p(z)^2\,dz}
\end{equation}
Using \(T_p(\cos\theta)=\cos(p\theta)\) and
\(T_p(\cosh t)=\cosh(pt)\) \cite[Eq.~18.5.1]{DLMF}, the denominator is
\(D_p\), and integration of \(\cosh^2(pt)\sinh t\) over \((0,a)\)
gives \(I_p\). Since
\(I_p\sim \sinh(a)\rho^{2p}/(8p)=2^{-3/2}\rho^{2p}/p\)
and \(D_p\to1\), this proves the lower bound and its asymptotic formula.

For the matching upper bound, write
\(\phi_j(z)=\sqrt{(2j+1)/2}\,P_j(z)\) for the orthonormal Legendre
basis on \((-1,1)\), and set \(\chi(z)=z+\sqrt{z^2-1}\), so that \(\chi(3)=\rho\).
The integral representation \cite[Eq.~14.12.7, \(m=0\)]{DLMF}
\begin{equation}
P_j(z)=\frac1\pi\int_0^\pi
\bigl(z+\sqrt{z^2-1}\cos\theta\bigr)^j\,d\theta
\end{equation}
implies \(|P_j(z)|\le\chi(2)^j\) on \([1,2]\), and
\(|P_j(z)|\le C\chi(z)^j/\sqrt{j+1}\) uniformly on \([2,3]\).
Indeed, put \(b(z)=\sqrt{z^2-1}/\chi(z)\ge b_0>0\) on \([2,3]\).
For \(j\ge1\), the inequalities \(1-t\le e^{-t}\) and
\(1-\cos\theta\ge2\theta^2/\pi^2\) give
\begin{equation}
\frac{|P_j(z)|}{\chi(z)^j}
\le\frac1\pi\int_0^\pi[1-b(z)(1-\cos\theta)]^j\,d\theta
\le\frac1\pi\int_0^\infty e^{-2jb_0\theta^2/\pi^2}\,d\theta
\le Cj^{-1/2}
\end{equation}
Increasing \(C\) covers \(j=0\).
Thus the reproducing kernel \(K_p(z,z)=\sum_{j=0}^p\phi_j(z)^2\)
satisfies the following bound. On \([1,2]\), use
\(\sum_{j=0}^p(2j+1)r^j\le(2p+1)r^p/(1-r^{-1})\), with
\(r=\chi(2)^2>1\), to retain only one power of \(p+1\):
\begin{equation}
K_p(z,z)\le
\begin{cases}
C(p+1)\chi(2)^{2p},&1\le z\le2\\
C\chi(z)^{2p},&2\le z\le3
\end{cases}
\end{equation}
Since \((\log\chi(z))'=1/\sqrt{z^2-1}\ge1/\sqrt8\) on \([2,3]\),
integration gives
\begin{equation}
\int_1^3K_p(z,z)\,dz\le
C(p+1)\chi(2)^{2p}+C\rho^{2p}/p
\le C\rho^{2p}/p
\end{equation}
Finally, Cauchy--Schwarz gives
\(|v(z)|^2\le K_p(z,z)\|v\|_{(-1,1)}^2\) for \(v\in\Pp_p\).
Integrating over \((1,3)\) proves the upper bound for \(C_p\).
\end{proof}

\begin{rem}
\label{rem:higher-dimensional-barrier}
The Chebyshev trial polynomial establishes the sharp exponential base and
algebraic prefactor, but need not maximize the norm quotient. In particular,
\eqref{eq:continuation-sharp-order} does not identify an asymptotic equivalent
for the optimal \(C_p\). The simpler choice \((x+1)^p\) still gives
\(C_p\ge\sqrt{2^{2p+1}-1}\), but loses the sharp exponential base.
The same constructive lower bounds apply on neighboring Cartesian cells in any
dimension by using a polynomial in the coordinate normal to the common face.
\end{rem}

To indicate where this constant enters, suppose a full-element extension estimate has
the schematic form
\begin{equation}
\|\nabla v_h\|_{\Omega_h}^2
\lesssim
C_p^2\|\nabla v_h\|_\Omega^2
+
\|v_h\|_{s_h}^2 \label{eq:full-element-extension}
\end{equation}
Combining \eqref{eq:full-element-extension} with the complete-element inverse trace
inequality gives
\begin{equation}
h\|\partial_n v_h\|_{\partial\Omega}^2
\lesssim
p^2\bigl(C_p^2\|\nabla v_h\|_\Omega^2+\|v_h\|_{s_h}^2\bigr) \label{eq:full-element-flux-bound}
\end{equation}
and the corresponding sufficient Nitsche penalty contains \(p^2C_p^2\). Proposition
\ref{prop:interval-continuation} shows that this proof mechanism yields an exponentially
large sufficient penalty.

This observation does not show that the discrete method itself requires an exponential
penalty. It shows that complete-element continuation is too strong an intermediate
estimate when the objective is only physical-boundary coercivity. The next section
replaces \eqref{eq:full-element-extension} by a one-sided estimate that never evaluates
the projected polynomial beyond the physical domain.

\section{One-Sided Stability Analysis}
\label{sec:one-sided-stability}

\subsection{Endpoint Inverse Control on Inward Segments}
\label{subsec:endpoint-inverse}

The \(p\)-dependence of the analysis is determined by a one-dimensional endpoint
estimate. We recall the classical Legendre and Jacobi endpoint argument
underlying polynomial inverse trace bounds; see
\cite{Schwab1998,WarburtonHesthaven2003}.

\begin{lem}[Endpoint Inverse Inequalities]
\label{lem:endpoint-inverse}
Let \(I_\delta=(0,\delta)\), \(\delta>0\), and \(r\in\Pp_m(I_\delta)\). Then
\begin{equation}
\delta |r(0)|^2
\le
(m+1)^2\|r\|_{I_\delta}^2 \label{eq:endpoint-inverse}
\end{equation}
More generally, if \(\alpha\) is a nonnegative integer and
\(\widehat r\in\Pp_m(0,1)\), then
\begin{equation}
|\widehat r(1)|^2
\le
(m+1)(m+\alpha+1)
\int_0^1 t^\alpha |\widehat r(t)|^2\,dt 
\label{eq:weighted-endpoint-inverse}
\end{equation}
Both estimates hold componentwise for vector-valued polynomials.
\end{lem}

\begin{proof}
Let \(P_j\) be the Legendre polynomial of degree \(j\) on \((-1,1)\). The functions
\begin{equation}
\varphi_j(t)
=
\left(\frac{2j+1}{\delta}\right)^{1/2}
P_j\left(\frac{2t}{\delta}-1\right),
\qquad
0\le j\le m \label{eq:shifted-legendre-basis}
\end{equation}
form an orthonormal basis of \(\Pp_m(I_\delta)\). Writing
\(r=\sum_{j=0}^m r_j\varphi_j\), the Cauchy--Schwarz inequality gives
\begin{align}
|r(0)|^2
&\le
\left(\sum_{j=0}^m r_j^2\right)
\left(\sum_{j=0}^m\varphi_j(0)^2\right) \label{eq:endpoint-cauchy}\\
&=
\frac{(m+1)^2}{\delta}\|r\|_{I_\delta}^2 \label{eq:endpoint-kernel}
\end{align}
where \(P_j(-1)^2=1\) and
\(\sum_{j=0}^m(2j+1)=(m+1)^2\) were used.

To prove \eqref{eq:weighted-endpoint-inverse}, let
\(P_j^{(0,\alpha)}\) be the Jacobi polynomial of degree \(j\). The shifted
polynomials
\begin{equation}
\psi_j(t)
=
(2j+\alpha+1)^{1/2}P_j^{(0,\alpha)}(2t-1),
\qquad 0\le j\le m
\label{eq:shifted-jacobi-basis}
\end{equation}
are orthonormal in \(L^2((0,1),t^\alpha\,dt)\), and
\(P_j^{(0,\alpha)}(1)=1\). Expanding \(\widehat r\) in this basis and applying the
Cauchy--Schwarz inequality therefore gives
\begin{align}
|\widehat r(1)|^2
&\le
\left(\int_0^1t^\alpha|\widehat r(t)|^2\,dt\right)
\sum_{j=0}^m\psi_j(1)^2 \\
&=
(m+1)(m+\alpha+1)
\int_0^1t^\alpha|\widehat r(t)|^2\,dt 
\end{align}
where
\(\sum_{j=0}^m(2j+\alpha+1)=(m+1)(m+\alpha+1)\).
\end{proof}

\begin{rem}[Sharp Endpoint Constant]
Equality in \eqref{eq:endpoint-inverse} is attained by the endpoint kernel
\(r(t)=\sum_{j=0}^m\varphi_j(0)\varphi_j(t)\).
Thus the squared constant \((m+1)^2\), and its square-root norm constant
\(m+1\), are sharp. The Jacobi endpoint kernel similarly attains
\eqref{eq:weighted-endpoint-inverse}.
\end{rem}

For the normal-flux estimate, we decompose \(v_h=q_{\mathcal P}+r_{\mathcal P}\),
with \(q_{\mathcal P}=E_{\mathcal P}v_h\) and
\(r_{\mathcal P}=v_h-q_{\mathcal P}\).
Figure~\ref{fig:one-sided-decomposition} illustrates the two controls, using
\(\nabla q_{\mathcal P}\) on \(U_{\mathcal P}\subset\Omega\) for the polynomial flux
and \(\nabla r_{\mathcal P}\) on the complete patch \(\mathcal P\) for the residual flux.
\begin{figure}[htbp]
\centering
\resizebox{\linewidth}{!}{\begingroup
  \fontsize{16}{19.2}\selectfont
  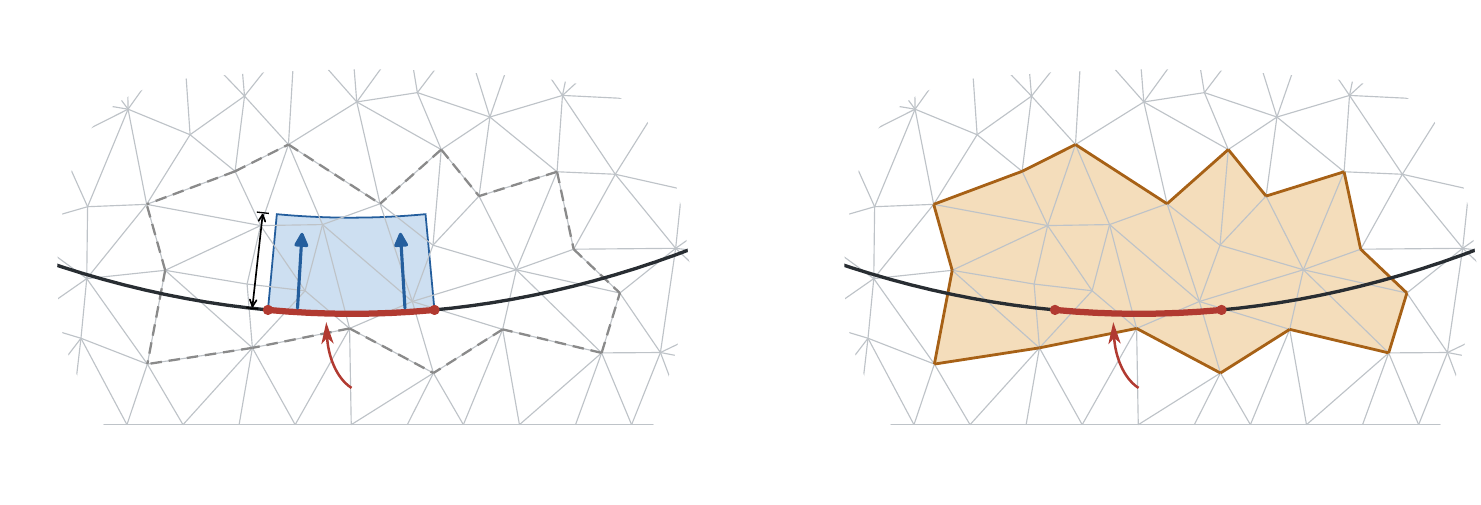\endgroup
}
\caption{One-sided flux control for \(v_h=q_{\mathcal P}+r_{\mathcal P}\).
The projected component \(q_{\mathcal P}=E_{\mathcal P}v_h\) is a single
polynomial on the whole patch, but its normal flux is estimated using its
gradient only in \(U_{\mathcal P}\subset\Omega\) (left).
One-dimensional endpoint estimates are applied along \(\Phi(y,t)\) for
\(y\in\Gamma_{\mathcal P}\) and \(0<t<c_\delta h\).
For \(r_{\mathcal P}\), complete-element inverse trace estimates use the
gradient over \(\mathcal P\), including exterior portions (right); this is
the gradient controlled by stabilization.
Both fluxes are estimated on the same boundary cell \(\Gamma_{\mathcal P}\).}
\label{fig:one-sided-decomposition}
\end{figure}

\begin{lem}[One-Sided Polynomial Flux Control]
\label{lem:polynomial-flux}
For each \(\mathcal P\in\Ph\) and each \(q\in\Pp_p(\mathcal P)\),
\begin{equation}
h\|\partial_n q\|_{\Gamma_{\mathcal P}}^2
\le
C_{\mathrm{line}}p^2\|\nabla q\|_{U_{\mathcal P}}^2,
\qquad
C_{\mathrm{line}}:=\frac{1}{c_\delta c_J}
\label{eq:polynomial-one-sided-flux}
\end{equation}
where \(c_\delta\) is the fixed tube-depth parameter and \(c_J\) is the lower
Jacobian bound in \eqref{eq:normal-coordinate-jacobian}.
\end{lem}

\begin{proof}
For \(y\in\Gamma_{\mathcal P}\), define
\begin{equation}
r_y(t)=q\bigl(\Phi(y,t)\bigr),
\qquad
0<t<c_\delta h \label{eq:line-polynomial}
\end{equation}
Because \(t\mapsto\Phi(y,t)\) is affine for fixed \(y\), \(r_y\) is a polynomial of
degree at most \(p\), and
\begin{equation}
r_y'(t)
=
-\normal(y)\cdot\nabla q\bigl(\Phi(y,t)\bigr) \label{eq:line-polynomial-derivative}
\end{equation}
Applying Lemma \ref{lem:endpoint-inverse} to \(r_y'\in\Pp_{p-1}\) gives
\begin{equation}
c_\delta h|\partial_nq(y)|^2
\le
p^2\int_0^{c_\delta h}
\bigl|\normal(y)\cdot\nabla q(\Phi(y,t))\bigr|^2\,dt \label{eq:pointwise-line-flux}
\end{equation}
Integrating over \(\Gamma_{\mathcal P}\), using the Jacobian bounds
\eqref{eq:normal-coordinate-jacobian}, and using
\(|\normal\cdot\nabla q|\le|\nabla q|\), yields
\eqref{eq:polynomial-one-sided-flux}.
\end{proof}

\subsection{Control of the Projection Residual}
\label{subsec:projection-residual}

We also use a cut-independent inverse trace inequality on a complete
background element. For curved cut sections, the \(p^2/h\) scaling follows
from the scalar inverse trace estimate of Wu and Xiao
\cite[Lemma 3.1(ii)]{WuXiao2019}: apply it to each component of the
polynomial gradient and enlarge the selected subcell to the complete element.
For simplex faces, sharp constants are given by Warburton and Hesthaven
\cite{WarburtonHesthaven2003}. The following formulation records an explicit
geometric constant for use in the penalty bound.

\begin{lem}[Complete-Element Inverse Trace]
\label{lem:complete-element-trace}
For \(0<h\le h_0\), every active element \(T\) with \(|\partial\Omega\cap T|>0\), and
\(w_h\in\Pp_p(T)\),
\begin{equation}
h\|\partial_n w_h\|_{\partial\Omega\cap T}^2
\le
C_{\mathrm{tr}}p^2\|\nabla w_h\|_T^2,
\qquad
C_{\mathrm{tr}}:=\frac{d}{c_*}=\frac{2d}{\varrho}
\label{eq:complete-element-inverse-trace}
\end{equation}
where every element contains a ball of radius \(\varrho h\),
\(c_*:=\varrho/2\), and \(h_0\) is chosen according to the boundary geometry
and mesh regularity. The constant is independent of the cut configuration,
\(h\), and \(p\).
\end{lem}

The proof with the stated geometric constant and mesh restriction is given
in Appendix~\ref{app:complete-element-trace}.

\subsection{The One-Sided Normal-Flux Estimate}
\label{subsec:normal-flux-estimate}

\begin{thm}[One-Sided Normal-Flux Inverse Estimate]
\label{thm:one-sided-flux}
For every \(0<h\le h_0\), the stabilization \eqref{eq:patch-stabilization}
satisfies
\begin{equation}
h\|\partial_n v_h\|_{\partial\Omega}^2
\le
C_{\mathrm{inv}}p^2
\bigl(\|\nabla v_h\|_\Omega^2+\|v_h\|_{s_h}^2\bigr)
\qquad \forall v_h\in V_h \label{eq:one-sided-flux}
\end{equation}
where one admissible constant is
\begin{equation}
C_{\mathrm{inv}}
:=
2C_{\mathrm{tr}}+4C_{\mathrm{line}}
\label{eq:assembled-inverse-constant}
\end{equation}
Thus \(C_{\mathrm{inv}}\) depends only on \(c_\delta,c_J,c_*\),
and \(d\), and is independent of the cut configuration, \(h\), and \(p\).
\end{thm}

\begin{proof}
Set \(r_{\mathcal P}=v_h-E_{\mathcal P}v_h\) and
\(q_{\mathcal P}=E_{\mathcal P}v_h\). On \(\Gamma_{\mathcal P}\), decompose
\begin{equation}
v_h
=r_{\mathcal P}+q_{\mathcal P} \label{eq:projection-decomposition}
\end{equation}
Figure~\ref{fig:one-sided-decomposition} shows the regions used to control these two components.
The patch also contains the complete boundary elements needed for the
residual trace: for almost every \(y\in\Gamma_{\mathcal P}\), the inward
points \(\Phi(y,t)\) lie, for sufficiently small \(t>0\), in an element
included in the patch construction. Its closure contains \(y\), with the
physical-side trace convention of Section~\ref{subsec:nitsche-formulation}.
Summing the complete-element inverse trace inequality over these elements gives
\begin{equation}
h\|\partial_n r_{\mathcal P}\|_{\Gamma_{\mathcal P}}^2
\le
C_{\mathrm{tr}}p^2\|\nabla r_{\mathcal P}\|_{\mathcal P}^2 
\label{eq:residual-flux-control}
\end{equation}
For the projected polynomial, Lemma \ref{lem:polynomial-flux} gives
\begin{align}
h\|\partial_nq_{\mathcal P}\|_{\Gamma_{\mathcal P}}^2
&\le
C_{\mathrm{line}}p^2\|\nabla q_{\mathcal P}\|_{U_{\mathcal P}}^2
\label{eq:projected-flux-line-control}\\
&\le
2C_{\mathrm{line}}p^2\left(
\|\nabla v_h\|_{U_{\mathcal P}}^2
+
\|\nabla r_{\mathcal P}\|_{\mathcal P}^2
\right) \label{eq:projected-flux-physical-control}
\end{align}
where \(q_{\mathcal P}=v_h-r_{\mathcal P}\) and
\(U_{\mathcal P}\subset\mathcal P\cap\Omega\) were used. Hence
\begin{align}
h\|\partial_n v_h\|_{\Gamma_{\mathcal P}}^2
&\le
p^2\Bigl(
4C_{\mathrm{line}}\|\nabla v_h\|_{U_{\mathcal P}}^2
+
(2C_{\mathrm{tr}}+4C_{\mathrm{line}})
\|\nabla r_{\mathcal P}\|_{\mathcal P}^2
\Bigr)
\label{eq:patch-flux-assembled}
\end{align}
The normal-coordinate map is a diffeomorphism and the boundary pieces
are disjoint up to surface-measure-zero sets. Their tubes are therefore
disjoint up to volume-measure-zero sets, so
\begin{equation}
\sum_{\mathcal P}\|\nabla v_h\|_{U_{\mathcal P}}^2
\le\|\nabla v_h\|_\Omega^2
\end{equation}
Summing over the boundary partition and recalling the definition of
\(\|v_h\|_{s_h}\) proves \eqref{eq:one-sided-flux} with
\eqref{eq:assembled-inverse-constant}.
\end{proof}

\begin{rem}[Alternative at Corners]
Physical cones can replace inward normal segments under a positive
height and bounded multiplicity assumption; see Appendix~\ref{app:cones}.
That variant retains an overlap factor unless a global multiplicity
bound is available.
\end{rem}

\begin{rem}[Boundary Assignment Leaves the Method Unchanged]
\label{rem:assignment-invariance}
The partition \(\{\Gamma_{\mathcal P}\}\) enters the proof, not the
bilinear form. To use an existing indexed patch family, partition the
boundary, up to surface-measure-zero sets, into measurable pieces \(A_j\)
contained in boundary elements \(T_j\). Choose a patch \(\sigma(j)\)
containing the complete element \(T_j\) and every inward segment of
depth \(c_\delta h\) based on \(A_j\). Set
\begin{equation}
\Gamma_{\mathcal P}=\bigcup_{\sigma(j)=\mathcal P}A_j
\end{equation}
With the same Jacobian, complete-element trace and patch-overlap bounds,
the preceding proof applies to this partition. Grouping the pieces by
patch before applying the trace estimate introduces no factor depending
on the number of pieces. Patches assigned an empty boundary piece may
retain their original stabilization terms. Thus reassignment changes
neither \(s_h\), the discrete system, the energy norm nor its consistency
error; it supplies the geometric decomposition used by the proof.
\end{rem}

\begin{rem}
\label{rem:no-full-element-control}
The full-element trace estimate is used only for the stabilized residual.
The projected polynomial is measured on \(U_{\mathcal P}\subset\Omega\),
so its continuation constant does not enter Theorem~\ref{thm:one-sided-flux}.
\end{rem}

\subsection{Coercivity and Continuity}
\label{subsec:coercivity-continuity}

Define the discrete energy norm
\begin{equation}
\energy{v}_h^2
=
\|\nabla v\|_\Omega^2
+
\|v\|_{s_h}^2
+
\frac{p^2}{h}\|v\|_{\partial\Omega}^2 \label{eq:energy-norm}
\end{equation}
and the augmented norm
\begin{equation}
\energy{v}_{h,*}^2
=
\|\nabla v\|_\Omega^2
+
\frac{p^2}{h}\|v\|_{\partial\Omega}^2
+
\frac{h}{p^2}\|\partial_n v\|_{\partial\Omega}^2 \label{eq:augmented-norm}
\end{equation}

\begin{thm}[Cut-Independent Coercivity]
\label{thm:coercivity}
For every \(0<h\le h_0\), a sufficient penalty condition is
\begin{equation}
\beta_0\ge\beta_*:=2C_{\mathrm{inv}}+\frac12 
\label{eq:penalty-threshold}
\end{equation}
Under this condition,
\begin{equation}
A_h(v_h,v_h)
\ge
\frac12\energy{v_h}_h^2
\qquad \forall v_h\in V_h \label{eq:coercivity}
\end{equation}
The constants are independent of \(h\), \(p\), and the cut configuration.
\end{thm}

\begin{proof}
Set
\(X=\|\nabla v_h\|_\Omega^2+\|v_h\|_{s_h}^2\) and
\(B=(p^2/h)\|v_h\|_{\partial\Omega}^2\). Theorem
\ref{thm:one-sided-flux} and Young's inequality give
\begin{align}
2|(\partial_n v_h,v_h)_{\partial\Omega}|
&\le
2(C_{\mathrm{inv}}XB)^{1/2} \label{eq:flux-product}\\
&\le
\frac12X+2C_{\mathrm{inv}}B  \label{eq:flux-young}
\end{align}
Substitution in \eqref{eq:stabilized-form} yields
\begin{equation}
A_h(v_h,v_h)
\ge
\frac12X+(\beta_0-2C_{\mathrm{inv}})B
\ge
\frac12(X+B)
\end{equation}
where the last inequality follows from \eqref{eq:penalty-threshold}. This is
\eqref{eq:coercivity}.
\end{proof}

\begin{rem}[Explicit Sufficient Penalty]
For \(h\le h_0\), the sufficient bound is
\begin{equation}
\beta_*=4C_{\mathrm{tr}}+8C_{\mathrm{line}}+\frac12
=\frac{8d}{\varrho}+\frac{8}{c_\delta c_J}+\frac12
\label{eq:explicit-sufficient-penalty}
\end{equation}
No patch-overlap factor enters this disjoint-tube estimate. Its slack comes
from the polynomial/residual splitting, the certified tube depth and the
complete-element trace constant.
\end{rem}

\begin{lem}[Nitsche Continuity]
\label{lem:nitsche-continuity}
For \(v\) with a square-integrable normal trace and \(w_h\in V_h\),
\begin{equation}
|a_h(v,w_h)|
\lesssim
\energy{v}_{h,*}\energy{w_h}_h \label{eq:nitsche-continuity}
\end{equation}
\end{lem}

\begin{proof}
Apply the Cauchy--Schwarz inequality to the volume and boundary-value terms. Scale the
term \((\partial_n v,w_h)_{\partial\Omega}\) with \(h^{1/2}/p\) and
\(p/h^{1/2}\). Scale \((v,\partial_n w_h)_{\partial\Omega}\) in the reverse order and
apply Theorem \ref{thm:one-sided-flux} to \(w_h\).
\end{proof}

\begin{rem}
\label{rem:condition-number}
The coercivity norm \eqref{eq:energy-norm} measures the physical gradient, the patch
residual, and the Nitsche trace. It is not uniformly equivalent to the Euclidean norm
of an arbitrary nodal coefficient vector. Consequently, Theorem \ref{thm:coercivity}
does not imply a cut-independent stiffness-matrix condition number. This limitation is
intentional and separates the present stability result from active-mesh extension and
preconditioning results.
\end{rem}

\section{\(hp\)-Error Analysis}
\label{sec:error-analysis}

\subsection{Approximation Tools}
\label{subsec:approximation-tools}

Throughout this section, fix \(\beta_0\ge\beta_*\); constants may depend
on this fixed penalty coefficient.
For the finite-regularity estimates, fix a bounded neighborhood \(\mathcal O\)
containing the background domain and all patch neighborhoods defined below.
Choose a fixed Sobolev extension operator \(E\) from \(\Omega\) to \(\mathcal O\).
For each fixed real \(s\ge2\), its boundedness gives
\begin{equation}
\|Eu\|_{H^s(\mathcal O)}
\le C_{E,s}\|u\|_{H^s(\Omega)}
\le C_{E,s}\|Eu\|_{H^s(\Omega_h)}
\label{eq:bounded-extension-control}
\end{equation}
The second inequality holds since \(\Omega\subset\Omega_h\) and \(Eu=u\)
on \(\Omega\). Here \(C_{E,s}\) is independent of \(h\), \(p\), and the cut position.
We use the Slobodeckij convention for noninteger Sobolev norms.

For the approximation tools, let \(s\ge2\) be an integer and set
\(\mu=\min(p+1,s)\). We use separate conforming and elementwise
approximants, as in \cite[Lemma 3.2]{WuXiao2019}, following
\cite{BabuskaSuri1987,Schwab1998}.
Conforming \(hp\)-approximation on the background mesh, restricted to the
active mesh, gives \(\pi_hEu\in V_h\) with
\begin{equation}
\|D^m(Eu-\pi_hEu)\|_{\Omega_h}
\lesssim
\frac{h^{\mu-m}}{p^{s-m}}\|Eu\|_{H^s(\Omega_h)}
\qquad m=0,1 \label{eq:hp-interpolation}
\end{equation}
The estimate on the fixed background domain implies this form by
\eqref{eq:bounded-extension-control}. Independently, choose a possibly
discontinuous piecewise polynomial \(\widetilde\pi_hEu\) satisfying
\begin{equation}
\|D^m(Eu-\widetilde\pi_hEu)\|_T
\lesssim
\frac{h^{\mu-m}}{p^{s-m}}\|Eu\|_{H^s(T)},
\qquad m=0,1,2,\quad T\in\Th
\label{eq:elementwise-hp-approximation}
\end{equation}
All second derivatives below are taken elementwise. The constants may depend
on the fixed index \(s\), including the extension norm \(C_{E,s}\), but not
on \(p\) or the cut configuration. We extend the final error estimates to
noninteger indices by interpolation on the fixed physical domain.

For \(\mathcal P=\mathcal P_i\), set
\(\mathcal N_{\mathcal P}=B(y_i,(C_0+1)h)\), with \(C_0\) from
\eqref{eq:patch-ball-containment}. Polynomial approximation on this ball,
followed by restriction to \(\mathcal P\), gives a polynomial
\(q_{\mathcal P}\in\Pp_p(\mathcal P)\) such that
\begin{equation}
\|\nabla(Eu-q_{\mathcal P})\|_{\mathcal P}
\lesssim
\frac{h^{\mu-1}}{p^{s-1}}
\|Eu\|_{H^s(\mathcal N_{\mathcal P})} \label{eq:patch-polynomial-approximation}
\end{equation}
The fixed radius factor is absorbed into the constant, since \(\mu\le s\).
Thus the constant may grow with the patch diameter relative to \(h\) and
with the fixed regularity index \(s\), but remains independent of \(p\).
The analytic estimate below treats the patch size explicitly instead of
substituting \(s=p+1\) into this constant.
The packing argument in Proposition~\ref{prop:boundary-patch-covering}
also bounds the overlap of these balls. Hence
\begin{equation}
\sum_{\mathcal P\in\Ph}\|Eu\|_{H^s(\mathcal N_{\mathcal P})}^2
\lesssim \|Eu\|_{H^s(\mathcal O)}^2
\lesssim \|Eu\|_{H^s(\Omega_h)}^2
\label{eq:patch-neighborhood-sum}
\end{equation}
The last inequality uses \eqref{eq:bounded-extension-control}; the balls
need not lie in \(\Omega_h\).

\begin{lem}[Approximation in the Augmented Norm]
\label{lem:augmented-approximation}
The interpolant in \eqref{eq:hp-interpolation} satisfies
\begin{equation}
\energy{Eu-\pi_hEu}_{h,*}
\lesssim
\frac{h^{\mu-1}}{p^{s-3/2}}
\|Eu\|_{H^s(\Omega_h)} \label{eq:augmented-approximation}
\end{equation}
\end{lem}

\begin{proof}
Set \(\eta_h=Eu-\pi_hEu\). The cut-element multiplicative trace inequality
\cite{HansboHansbo2002} and
\eqref{eq:hp-interpolation} give
\begin{align}
\|\eta_h\|_{\partial\Omega}^2
&\lesssim
h^{-1}\|\eta_h\|_{\Omega_h}^2
+
\|\eta_h\|_{\Omega_h}\|\nabla\eta_h\|_{\Omega_h} \label{eq:value-trace-interpolation}\\
&\lesssim
\frac{h^{2\mu-1}}{p^{2s-1}}
\|Eu\|_{H^s(\Omega_h)}^2 \label{eq:value-trace-rate}
\end{align}
For the normal trace, write
\(\eta_h=e_h^{\rm loc}+d_h\), where
\(e_h^{\rm loc}=Eu-\widetilde\pi_hEu\) and
\(d_h=\widetilde\pi_hEu-\pi_hEu\).
Apply the multiplicative trace inequality to \(\nabla e_h^{\rm loc}\)
and Lemma~\ref{lem:complete-element-trace} to the polynomial \(d_h|_T\).
Summing over cut elements gives
\begin{align}
\|\partial_n\eta_h\|_{\partial\Omega}^2
&\lesssim
h^{-1}\|\nabla e_h^{\rm loc}\|_{\Omega_h}^2
+
\|\nabla e_h^{\rm loc}\|_{\Omega_h}\|D^2e_h^{\rm loc}\|_{\Omega_h}
+\frac{p^2}{h}\|\nabla d_h\|_{\Omega_h}^2
\label{eq:normal-trace-interpolation}\\
&\lesssim
\frac{h^{2\mu-3}}{p^{2s-4}}
\|Eu\|_{H^s(\Omega_h)}^2 \label{eq:normal-trace-rate}
\end{align}
Here \eqref{eq:hp-interpolation} and
\eqref{eq:elementwise-hp-approximation} bound \(\nabla d_h\) by the sum
of the two gradient approximation errors. After multiplication by
\(h/p^2\), the normal-trace term is bounded by
\(h^{2\mu-2}p^{-2s+2}\|Eu\|_{H^s(\Omega_h)}^2\), which is of higher
order in \(p\) than the boundary-value contribution.
Combining \eqref{eq:hp-interpolation}, \eqref{eq:value-trace-rate}, and
\eqref{eq:normal-trace-rate} with \eqref{eq:augmented-norm} proves
\eqref{eq:augmented-approximation}.
\end{proof}

\begin{lem}[Approximation in the Stabilization Seminorm]
\label{lem:stabilization-approximation}
The exact extension and the interpolant in \eqref{eq:hp-interpolation} satisfy
\begin{equation}
\|Eu\|_{s_h}
+
\|\pi_hEu\|_{s_h}
+
\|Eu-\pi_hEu\|_{s_h}
\lesssim
\frac{h^{\mu-1}}{p^{s-1}}
\|Eu\|_{H^s(\Omega_h)} \label{eq:stabilization-approximation}
\end{equation}
\end{lem}

\begin{proof}
The \(H^1\)-projection property and polynomial reproduction give
\begin{align}
\|\nabla(\pi_hEu-E_{\mathcal P}\pi_hEu)\|_{\mathcal P}
&\le
\|\nabla(\pi_hEu-q_{\mathcal P})\|_{\mathcal P} \label{eq:projection-best-approximation}\\
&\le
\|\nabla(\pi_hEu-Eu)\|_{\mathcal P}
+
\|\nabla(Eu-q_{\mathcal P})\|_{\mathcal P} \label{eq:projection-error-split}
\end{align}
Use \eqref{eq:hp-interpolation} and
\eqref{eq:patch-polynomial-approximation}, square, sum over patches, and invoke
bounded patch overlap and \eqref{eq:patch-neighborhood-sum} to bound
\(\|\pi_hEu\|_{s_h}\). The same best-approximation argument applied
directly to \(Eu\) bounds \(\|Eu\|_{s_h}\). The triangle inequality then bounds
\(\|Eu-\pi_hEu\|_{s_h}\) and proves \eqref{eq:stabilization-approximation}.
\end{proof}

\subsection{Strang Estimate and Energy Error}
\label{subsec:strang-energy}

\begin{thm}[Best Approximation with Weak Consistency]
\label{thm:strang}
Let \(u\in H^{3/2+\varepsilon}(\Omega)\), \(\varepsilon>0\), solve
\eqref{eq:model-bulk}--\eqref{eq:model-boundary}, let \(Eu\) be a fixed
extension with the same regularity, and let \(u_h\) solve
\eqref{eq:discrete-method} with \(\beta_0\ge\beta_*\). Then
\begin{equation}
\energy{Eu-u_h}_h
\lesssim
\inf_{v_h\in V_h}
\left(
\energy{Eu-v_h}_{h,*}
+
\|Eu-v_h\|_{s_h}
+
\|v_h\|_{s_h}
\right) \label{eq:strang-estimate}
\end{equation}
\end{thm}

\begin{proof}
For arbitrary \(v_h\in V_h\), set \(e_h=v_h-u_h\). Consistency of the Nitsche
form and the discrete equation give
\begin{equation}
A_h(e_h,e_h)
=
a_h(v_h-Eu,e_h)+s_h(v_h,e_h) \label{eq:strang-identity}
\end{equation}
Coercivity, Lemma \ref{lem:nitsche-continuity}, and the Cauchy--Schwarz inequality in
the stabilization seminorm imply
\begin{equation}
\energy{e_h}_h
\lesssim
\energy{Eu-v_h}_{h,*}+\|v_h\|_{s_h} \label{eq:discrete-error-bound}
\end{equation}
Moreover,
\begin{equation}
\energy{Eu-v_h}_h
\le
\energy{Eu-v_h}_{h,*}+\|Eu-v_h\|_{s_h} \label{eq:energy-augmented-comparison}
\end{equation}
The triangle inequality, \eqref{eq:discrete-error-bound}, and minimization over
\(v_h\) prove \eqref{eq:strang-estimate}.
\end{proof}

\begin{thm}[\(hp\)-Energy Estimate]
\label{thm:hp-energy}
Let \(Eu\) be the fixed bounded extension above, with
\(Eu\in H^s(\Omega_h)\) for a real \(s\ge2\), and let
\(\mu=\min(p+1,s)\). Then
\begin{equation}
\energy{Eu-u_h}_h
\lesssim
\frac{h^{\mu-1}}{p^{s-3/2}}
\|Eu\|_{H^s(\Omega_h)} \label{eq:hp-energy-error}
\end{equation}
with a constant that may depend on \(s\), but is independent of the cut configuration, \(h\), and \(p\).
\end{thm}

\begin{proof}
First let \(s\) be an integer. Choose \(v_h=\pi_hEu\) in
Theorem~\ref{thm:strang}. Apply Lemmas
\ref{lem:augmented-approximation} and \ref{lem:stabilization-approximation}, and use the bounded overlap and patch diameter bounds in Proposition
\ref{prop:boundary-patch-covering}.

For \(s=k+\theta\), \(k\ge2\) integer and \(0<\theta<1\), consider the
linear error operator \(u\mapsto Eu-u_h(u)\), where \(u_h(u)\) uses the
data \(f=-\Delta u\) and \(g=u|_{\partial\Omega}\). For fixed \(h,p\),
the energy norm defines the same target space at both regularity indices.
By \eqref{eq:bounded-extension-control}, the integer estimates at \(k\) and
\(k+1\) hold with \(\|u\|_{H^k(\Omega)}\) and
\(\|u\|_{H^{k+1}(\Omega)}\) on the right. Real interpolation
\((\cdot,\cdot)_{\theta,2}\), using
\((H^k(\Omega),H^{k+1}(\Omega))_{\theta,2}=H^{k+\theta}(\Omega)\)
with equivalent norms on the fixed Lipschitz domain \(\Omega\)
\cite[Corollary 4.7]{ChandlerWildeHewettMoiola2015}, therefore gives the bound with
\(\|u\|_{H^s(\Omega)}\le\|Eu\|_{H^s(\Omega_h)}\).
The exponents interpolate to those in \eqref{eq:hp-energy-error}, because
the only breakpoint of \(\min(p+1,s)\) is the integer \(p+1\).
The interpolation constants depend on \(s\) and \(\Omega\), not on the cut.
\end{proof}

For analytic data we use a direct approximation construction, rather
than substituting \(s=p+1\) in an estimate whose constant depends on \(s\).
Only here it is useful to retain the patch scale
\begin{equation}
\widetilde h=\max_{\mathcal P\in\Ph}\operatorname{diam}(\mathcal P)
\simeq h
\label{eq:patch-mesh-scale}
\end{equation}
because a single polynomial must approximate the extension across each patch.

\begin{cor}[Analytic Data]
\label{cor:analytic-convergence}
Assume that the fixed extension is analytic in a neighborhood containing
all element stars and patch neighborhoods, with
\begin{equation}
\|D^\alpha Eu\|_{L^\infty}
\le C_u A^{|\alpha|}|\alpha|!
\qquad\text{for every multi-index }\alpha
\label{eq:analytic-regularity}
\end{equation}
There is a geometry constant \(C_a\), independent of the cut, \(h\) and
\(p\), such that, if \(\theta=C_a A\widetilde h\le\theta_0<1\), then
\begin{equation}
\energy{Eu-u_h}_h
\lesssim C_u h^{-1}p^{3/2}\theta^p 
\label{eq:analytic-energy-error}
\end{equation}
Here \(C_a\) absorbs the dimension-dependent conversion between
coordinate and Euclidean distances, element-star and patch diameters
relative to \(\widetilde h\), and mesh shape regularity and quasi-uniformity.
The constant may depend on \(\theta_0\) but not on \(p\). Consequently
the error decays exponentially with \(p\) on each such fixed mesh.
\end{cor}

\begin{proof}
Let \(\lambda_a\) be the piecewise affine vertex functions on the
background mesh and let \(T_{a,p-1}Eu\) be the Taylor polynomial at vertex
\(a\) of degree \(p-1\). The function
\begin{equation}
v_h=\sum_a\lambda_a T_{a,p-1}Eu
\end{equation}
is continuous and piecewise of degree at most \(p\).
The Taylor series, \eqref{eq:analytic-regularity}, the bounded number of
vertices per element and \(\|\nabla\lambda_a\|_\infty\lesssim h^{-1}\)
give, for \(m=0,1,2\), the broken-norm estimates
\begin{equation}
\|D^m(Eu-v_h)\|_{\Omega_h}
\lesssim C_u h^{-m}p^m\theta^p
\end{equation}
Indeed, differentiating a Taylor tail at most twice introduces at most a
quadratic factor in its order; summing the resulting geometric series is
uniform for \(\theta\le\theta_0<1\). Enlarging \(C_a\) covers the
element stars and the finitely overlapping patch neighborhoods.
Applying the cut-element multiplicative trace inequality to \(Eu-v_h\)
and its elementwise gradient therefore yields
\(\energy{Eu-v_h}_{h,*}\lesssim C_u h^{-1}p^{3/2}\theta^p\).
A Taylor polynomial of degree \(p-1\) on each patch gives the same
gradient approximation with bound \(C_u h^{-1}p\theta^p\).
Projection best approximation, the triangle inequality and bounded overlap
then bound \(\|Eu-v_h\|_{s_h}+\|v_h\|_{s_h}\) by that quantity.
Theorem~\ref{thm:strang} proves the displayed bound. No smoothness-dependent interpolation constant is evaluated at
\(s=p+1\).
\end{proof}

\subsection{\(L^2\)-Error Estimate}
\label{subsec:l2-error}

Because \(\Omega\) is \(C^2\), the dual Dirichlet problem has \(H^2\)-regularity.
Let \(z\in H^2(\Omega)\cap H_0^1(\Omega)\) solve
\begin{equation}
-\Delta z=u-u_h
\quad\text{in }\Omega \label{eq:dual-problem}
\end{equation}
and let \(Ez\) be a fixed \(H^2\)-extension. Then
\begin{equation}
\|Ez\|_{H^2(\Omega_h)}
\lesssim
\|u-u_h\|_\Omega \label{eq:dual-regularity}
\end{equation}

\begin{thm}[\(L^2\)-Error Estimate]
\label{thm:l2-error}
Under the assumptions of Theorem~\ref{thm:hp-energy}, let the fixed bounded
extension satisfy \(Eu\in H^s(\Omega_h)\) for a real \(s\ge2\), and set
\(\mu=\min(p+1,s)\). With the dual regularity above,
\begin{equation}
\|u-u_h\|_\Omega
\lesssim
\frac{h}{p^{1/2}}
\frac{h^{\mu-1}}{p^{s-3/2}}
\|Eu\|_{H^s(\Omega_h)} \label{eq:l2-error}
\end{equation}
and therefore
\begin{equation}
\|u-u_h\|_\Omega
\lesssim
\frac{h^\mu}{p^{s-1}}
\|Eu\|_{H^s(\Omega_h)} \label{eq:l2-error-simplified}
\end{equation}
\end{thm}

\begin{proof}
First let \(s\) be an integer. Set \(z_h=\pi_hEz\) and
\(\xi_h=\pi_hEu-u_h\). Integration by parts in the dual
problem and the discrete equation give
\begin{align}
\|u-u_h\|_\Omega^2
&=
a_h(u-u_h,z) \label{eq:dual-identity}\\
&=
a_h(u-u_h,z-z_h)+s_h(u_h,z_h) \label{eq:dual-galerkin-identity}
\end{align}
The first term on the right is decomposed as
\begin{equation}
a_h(u-u_h,z-z_h)
=
a_h(Eu-\pi_hEu,Ez-z_h)+a_h(\xi_h,Ez-z_h) \label{eq:dual-nitsche-split}
\end{equation}
The Nitsche form is continuous in the product of two augmented norms, and its
symmetry together with Lemma \ref{lem:nitsche-continuity} gives
\begin{align}
|a_h(u-u_h,z-z_h)|
&\lesssim
\left(
\energy{Eu-\pi_hEu}_{h,*}
+
\energy{\xi_h}_h
\right)
\energy{Ez-z_h}_{h,*} \label{eq:dual-nitsche-bound}\\
|s_h(u_h,z_h)|
&\le
\left(
\|\pi_hEu\|_{s_h}
+
\|\xi_h\|_{s_h}
\right)
\|z_h\|_{s_h} \label{eq:dual-stabilization-bound}
\end{align}
Lemma \ref{lem:augmented-approximation} with \(s=2\), Lemma
\ref{lem:stabilization-approximation}, and dual regularity yield
\begin{equation}
\|z_h\|_{s_h}
+
\energy{Ez-z_h}_{h,*}
\lesssim
\frac{h}{p^{1/2}}\|u-u_h\|_\Omega \label{eq:dual-approximation}
\end{equation}
Theorem \ref{thm:strang} and the triangle inequality bound
\(\energy{\xi_h}_h\) and \(\|\xi_h\|_{s_h}\) by the primal approximation terms.
Consequently, \eqref{eq:dual-galerkin-identity}--\eqref{eq:dual-approximation}
give
\begin{align}
\|u-u_h\|_\Omega^2
&\lesssim
\frac{h}{p^{1/2}}
\left(
\energy{Eu-\pi_hEu}_{h,*}
+
\|\pi_hEu\|_{s_h}
\right)
\|u-u_h\|_\Omega \label{eq:dual-error-bound}\\
&\lesssim
\frac{h}{p^{1/2}}
\frac{h^{\mu-1}}{p^{s-3/2}}
\|Eu\|_{H^s(\Omega_h)}
\|u-u_h\|_\Omega \label{eq:dual-final-bound}
\end{align}
Division by \(\|u-u_h\|_\Omega\) proves \eqref{eq:l2-error}; the simplified form
follows by collecting the powers of \(h\) and \(p\).
For noninteger \(s\), interpolate the linear map
\(u\mapsto u-u_h(u)\) into \(L^2(\Omega)\) between neighboring integer
indices, exactly as in the proof of Theorem~\ref{thm:hp-energy}.
\end{proof}

\begin{rem}[Degree Suboptimality]
Compared with the optimal approximation orders
\(h^{\mu-1}/p^{s-1}\) in the gradient norm and
\(h^\mu/p^s\) in \(L^2\), the energy estimate loses half a power of
\(p\) and \eqref{eq:l2-error-simplified} loses a full power. The latter
combines the half-power loss in primal approximation with another in the
dual estimate. These are upper bounds; faster observed rates do not
establish that either loss is necessary. A corresponding half-power
degree loss occurs in classical interior-penalty \(hp\)-DG estimates
\cite{HoustonSchwabSuli2002} and in the broken \(H^1\)-error estimate of
Chernov and Hansbo~\cite[Theorem 3]{ChernovHansbo2011} for an \(hp\)-Nitsche
method on interface-fitted, nonmatching subdomain meshes.
\end{rem}

\section{Numerical Results}
\label{sec:numerical-results}
The numerical companion~\cite{HpNumericalCompanion2026} contains the
implementation, complete result tables, and verification records for all
experiments in this section. Appendix~\ref{sec:numerical-details} summarizes
the verification procedures and provides a guide to the supporting material.

All experiments are two-dimensional. The domain is the disk
\(\Omega_{c}=\{x:|x-c|<0.62\}\subset(-1,1)^2\),
and the background square is divided into \(2N^2\) affine triangles,
with \(h=2/N\). The conforming hierarchical \(\Pp_p\) basis consists of
barycentric vertex functions, consistently oriented integrated-Legendre
edge modes, and interior bubbles. The vertex and edge construction follows
the hierarchical \(H^1\) family in
\cite[Section II]{SchoeberlZaglmayr2005}. The bubbles span
\(\lambda_1\lambda_2\lambda_3\Pp_{p-3}\), where \(\lambda_i\) are barycentric
coordinates, and are orthogonalized within that space in the reference
\(H^1\) inner product. This invertible basis change preserves \(\Pp_p\);
for \(p<3\), the bubble space is empty. The vertex and edge functions span
the polynomial boundary traces, while the bubbles span the zero-trace
remainder; together they span \(\Pp_p\). The guide details the construction,
spanning argument and implementation checks.

Boundary patches contain complete active triangles within \(L\) face-adjacency steps;
\(L\) denotes the patch layer count, while \(L_{\rm pc}\) below denotes a
preconditioner block depth. Unless specified, \(L=2\), the stabilization
is rescaled by \(\gamma_1=0.1\), and \(\beta_0=16\). This positive
rescaling changes constants, not the degree dependence of the theory:
for $\gamma>0$, $C_\gamma^{\rm num}\le\max(1,\gamma^{-1})C_1^{\rm num}$,
since $K+\gamma S\succeq\min(1,\gamma)(K+S)$.
The fixed centers are \(c_A=(0.071,-0.043)\) and
\(c_B=(-0.113,0.087)\). No known-solution stabilization load is added.

The baseline PDE computations use double-precision assembly with stable
patch projections and iterative refinement using extended-precision residuals.
The high-accuracy points in Figure~\ref{fig:sequence-p} additionally use
higher-precision analytic assembly and retain extra precision in the solution
coefficients to resolve the smallest errors.

\subsection{Continuation and the One-Sided Flux Constant}
\label{subsec:numerical-continuation}
The computed continuation constant is the square root of the largest
generalized eigenvalue of the two polynomial Gram matrices on
\((-1,0)\) and \((0,1)\). High-precision repeats confirm
\(C_{14}=9.57615\times10^9\), compared with the Chebyshev lower bound
\(B_{14}=8.14801\times10^9\). Figure~\ref{fig:continuation-constants}
shows the sharp exponential base alongside the linear endpoint norm
constant \(p+1\); at \(p=14\), the optimum exceeds the Chebyshev lower
bound by about 18\%, without asserting an asymptotic ratio.

\begin{figure}[!htbp]\centering
\begin{tikzpicture}\begin{axis}[width=.47\textwidth,height=.27\textwidth,font=\small,tick label style={font=\scriptsize},legend style={font=\scriptsize,draw=none,fill=none,at={(0.5,1.06)},anchor=south,/tikz/column sep=3pt},legend columns=-1,grid=major,ymode=log,xlabel={$p$},ylabel={norm constant}]
\addplot[blue,mark=*] coordinates {(1,3.73205080757) (2,15.9372539332) (3,76.6475588297) (4,389.480940998) (5,2039.37433456) (6,10883.9749885) (7,58862.3277644) (8,321464.201888) (9,1768842.24017) (10,9791075.56009) (11,54459188.1503) (12,304124359.731) (13,1704102991.36) (14,9576145559.79)};
\addlegendentry{$C_p$}
\addplot[red,dashed,mark=none] coordinates {(1,3.60555127546) (2,13.1312060158) (3,63.4669900305) (4,324.648445096) (5,1708.20156651) (6,9149.60009208) (7,49620.0606104) (8,271582.831149) (9,1497011.51748) (10,8298452.07351) (11,46213147.862) (12,258341073.908) (13,1448848075.07) (14,8148014410.02)};
\addlegendentry{$B_p$}
\addplot[black,mark=square*] coordinates {(1,2) (2,3) (3,4) (4,5) (5,6) (6,7) (7,8) (8,9) (9,10) (10,11) (11,12) (12,13) (13,14) (14,15)};
\addlegendentry{$p+1$}
\end{axis}\end{tikzpicture}
\caption{Continuation grows with asymptotic slope $\log_{10}(3+2\sqrt2)\simeq0.766$ per degree, apart from the algebraic prefactor; endpoint control grows linearly. The $B_p$ and $C_p$ curves nearly overlap on this logarithmic scale; $C_{14}/B_{14}=1.175$.}
\label{fig:continuation-constants}
\end{figure}
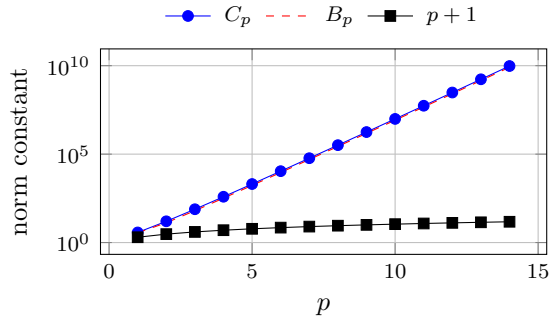

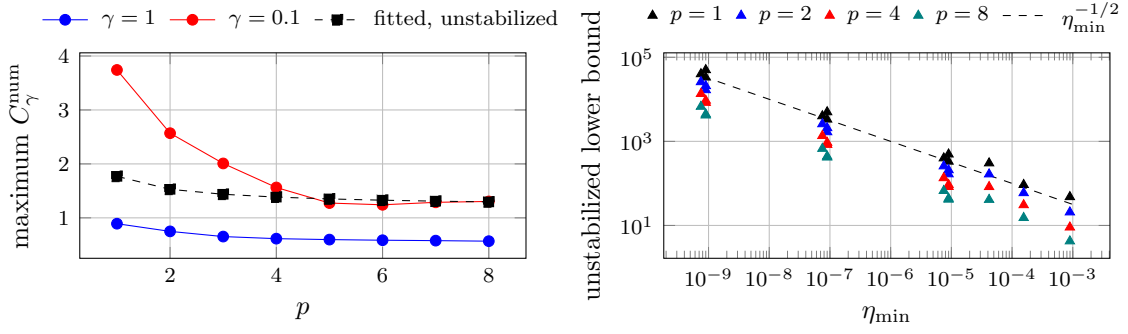
\begin{figure}[!tb]\centering
\begin{tikzpicture}[baseline=(current axis.south)]\begin{axis}[width=.47\textwidth,height=.27\textwidth,font=\small,tick label style={font=\scriptsize},legend style={font=\scriptsize,draw=none,fill=none,at={(0.5,1.06)},anchor=south,/tikz/column sep=3pt},legend columns=-1,grid=major,xlabel={$p$},ylabel={maximum $C_\gamma^{\rm num}$}]
\addplot[blue,mark=*] coordinates {(1,0.891617350655) (2,0.749331927624) (3,0.65309554804) (4,0.614496513503) (5,0.596974356042) (6,0.585366643541) (7,0.577728563003) (8,0.56780049005)};
\addlegendentry{$\gamma=1$}
\addplot[red,mark=*] coordinates {(1,3.74138291388) (2,2.56944474279) (3,2.00833338682) (4,1.56286935224) (5,1.27424756822) (6,1.24200139526) (7,1.2878421337) (8,1.30423670079)};
\addlegendentry{$\gamma=0.1$}
\addplot[black,dashed,mark=square*] coordinates {(1,1.76958499238) (2,1.52907475814) (3,1.43715313309) (4,1.38462057444) (5,1.35083646428) (6,1.32716583324) (7,1.30953415255) (8,1.29580720161)};
\addlegendentry{fitted, unstabilized}
\end{axis}\end{tikzpicture}\begin{tikzpicture}[baseline=(current axis.south)]\begin{axis}[width=.47\textwidth,height=.27\textwidth,font=\small,tick label style={font=\scriptsize},legend style={font=\scriptsize,draw=none,fill=none,at={(0.5,1.06)},anchor=south,/tikz/column sep=3pt},legend columns=-1,grid=major,xmode=log,ymode=log,xlabel={$\eta_{\min}$},ylabel={unstabilized lower bound}]
\addplot[black,only marks,mark=triangle*] coordinates {(0.000895218714972,47.645958703) (8.99516713225e-06,492.831779744) (8.99951616512e-08,4977.49523874) (8.99995220058e-10,49928.9554756) (0.000155285631915,92.6312908983) (8.99516713253e-06,332.706982595) (8.99951615729e-08,3332.70618229) (8.99994828663e-10,33332.7061014) (4.18446220432e-05,302.891655684) (7.49921615693e-06,398.636232671) (7.49992163998e-08,3998.62975772) (7.49999062855e-10,39998.6291047)};
\addlegendentry{$p=1$}
\addplot[blue,only marks,mark=triangle*] coordinates {(0.000895218714972,20.6389325266) (8.99516713225e-06,208.136762362) (8.99951616512e-08,2083.1365449) (8.99995220058e-10,20833.1365235) (0.000155285631915,59.2824157935) (8.99516713253e-06,166.264284344) (8.99951615729e-08,1666.26352538) (8.99994828663e-10,16666.2634495) (4.18446220432e-05,165.459151274) (7.49921615693e-06,257.083745087) (7.49992163998e-08,2582.07632505) (7.49999062855e-10,25832.075577)};
\addlegendentry{$p=2$}
\addplot[red,only marks,mark=triangle*] coordinates {(0.000895218714972,9.02334692059) (8.99516713225e-06,92.7701470811) (8.99951616512e-08,934.508210575) (8.99995220058e-10,9365.63765558) (0.000155285631915,30.6741044551) (8.99516713253e-06,83.0437551808) (8.99951615729e-08,833.043180376) (8.99994828663e-10,8333.04312354) (4.18446220432e-05,82.7385382929) (7.49921615693e-06,134.550848081) (7.49992163998e-08,1352.47727602) (7.49999062855e-10,13531.7907058)};
\addlegendentry{$p=4$}
\addplot[teal,only marks,mark=triangle*] coordinates {(0.000895218714972,4.22586779385) (8.99516713225e-06,43.5429283122) (8.99951616512e-08,440.282033999) (8.99995220058e-10,4419.28420593) (0.000155285631915,15.1870349426) (8.99516713253e-06,41.5108372625) (8.99951615729e-08,416.510532183) (8.99994828663e-10,4166.51050317) (4.18446220432e-05,41.0159642676) (7.49921615693e-06,66.6976376105) (7.49992163998e-08,670.66101039) (7.49999062855e-10,6710.32208158)};
\addlegendentry{$p=8$}
\addplot[black,dashed] coordinates {(8e-10,35355.3390593) (0.001,31.6227766017)};
\addlegendentry{$\eta_{\min}^{-1/2}$}
\end{axis}\end{tikzpicture}
\caption{Left: maximum over twelve cuts, and the unstabilized fitted-polygon comparator. Right: individual-basis Rayleigh lower-bound witnesses without stabilization; all three paths and an $\eta_{\min}^{-1/2}$ reference line are shown. These finite measurements test the normalization $h/p^2$ directly.}
\label{fig:direct-flux}
\end{figure}

We next measure Theorem~\ref{thm:one-sided-flux} directly. With
\(H_{ij}=(\partial_n\phi_i,\partial_n\phi_j)_{\partial\Omega}\),
\(K_{ij}=(\nabla\phi_i,\nabla\phi_j)_\Omega\) and raw patch matrix \(S\),
we compute
\begin{equation}
C_\gamma^{\rm num}=\lambda_{\max}\!\left(\frac{h}{p^2}H,K+\gamma S\right)
\quad\hbox{on }V_h/\mathbb R
\end{equation}
The derivative Gram matrix \(H\) is assembled separately from the boundary
value mass matrix. We test \(p=1,\ldots,8\) on all twelve disk placements
defined below, with \(\gamma=1,0.1,0\). Constants are removed by pinning
one vertex coefficient. Both quadratic forms are invariant under adding a
constant, so this selects one representative of each class in $V_h/\mathbb R$;
no shift is added. Two quadrature rules and
eigenpair checks distinguish resolved eigenvalues from lower bounds.
The unit and computational stabilizations give values in
\([0.560,0.892]\) and \([0.735,3.741]\), respectively.
Their maxima remain bounded over the tested degree range
(Figure~\ref{fig:direct-flux}); the theorem requires boundedness, not a
flat curve. For \(\gamma=0\), even individual basis-function Rayleigh
quotients grow strongly as cells become small. Where inversion of the
unstabilized metric is unresolved, these positive-integral lower bounds
are reported instead of a regularized eigenvalue. Here ``positive-integral
lower bounds'' means the maximum of the Rayleigh quotients of individual
basis functions, with numerator and denominator evaluated as positive
quadrature integrals; these are numerical lower-bound witnesses, not
interval-arithmetic enclosures. The maximization includes all basis functions.
The sliver scaling is illustrated by an exterior-vertex hat function with
physical support area \(\asymp\epsilon^2h^2\), boundary length
\(\asymp\epsilon h\), and normal derivative of size \(h^{-1}\) on a
comparable boundary portion. Its quotient is
\(\asymp\epsilon^{-1}\asymp\eta_{\min}^{-1/2}\) at fixed degree.
An interior-vertex hat need not have this small support.
The dashed guide in Figure~\ref{fig:direct-flux} shows this scaling.

The measured unit-stabilization maximum gives
$2C_1^{\rm num}+1/2\le2.284$ in the argument of
Theorem~\ref{thm:coercivity}, a measured sufficient penalty for the sampled
spaces. By contrast, \eqref{eq:explicit-sufficient-penalty} gives the
analytic unit-stabilization bound \(4C_{\mathrm{tr}}+256.5\) for the certified
disk tube parameters. Its explicit complete-element trace constant requires
the sufficient mesh condition \(h\le\varrho R/6\), which none of the sampled
meshes meets; Appendix~\ref{sec:numerical-details} gives the calculation.
Thus this analytic bound does not certify the reported thresholds.
The measurements below use the rescaled computational stabilization and
are not evaluations of that sufficient bound or guarantees for untested cuts.

For the computational stabilization, write
$X=v^T(K+0.1S)v$ and $B=(p^2/h)\|v\|_{\partial\Omega}^2$.
Cauchy--Schwarz gives
\begin{equation}
A(\beta_0)(v,v)\ge X-2(C_{0.1}^{\rm num}XB)^{1/2}+\beta_0B
=(\sqrt X-\sqrt{C_{0.1}^{\rm num}B})^2
+(\beta_0-C_{0.1}^{\rm num})B
\end{equation}
Thus the exact definiteness threshold satisfies
$\beta_{0,*}\le C_{0.1}^{\rm num}$; strict positivity follows for
$\beta_0>\nobreak C_{0.1}^{\rm num}$, including the nonzero constant mode.
All 36 matched cut measurements respect this ordering, even for the
reported upper endpoints with their $10^{-8}$ crossing margin. For example,
the oblique $p=2$, $\epsilon=3\times10^{-3}$ case has upper endpoint
$1.9132$ and $C_{0.1}^{\rm num}=2.5694$.
Figures~\ref{fig:direct-flux} and~\ref{fig:sequence-cut} therefore test the
same flux-control mechanism at two levels.

\subsection{Penalty Thresholds Across Cuts and Degrees}
\label{subsec:numerical-cut-stability}
At \(N=6\), translate the disk along directions proportional to
\((1,0),(1,1),(1,2)\), called axis, diagonal and oblique. A selected mesh
vertex has distance \(R-\epsilon h\) from the center, for
\(\epsilon\in\{3\times10^{-2},3\times10^{-3},3\times10^{-4},3\times10^{-5}\}\).
Cut severity is \(\eta_{\min}=\min_T|T\cap\Omega|/|T|\) over active cells.
The threshold study uses \(p=2,4,6\) and \(L=2\), with \(L=4\) controls
at the two endpoint offsets for \(p=4,6\). Geometry is fixed across degrees.

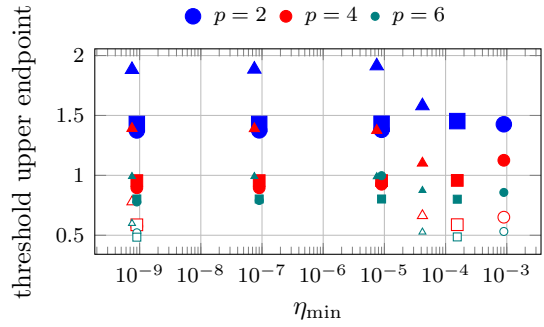
\begin{figure}[!htbp]\centering
\begin{tikzpicture}\begin{axis}[width=.47\textwidth,height=.27\textwidth,font=\small,tick label style={font=\scriptsize},legend style={font=\scriptsize,draw=none,fill=none,at={(0.5,1.06)},anchor=south,/tikz/column sep=3pt},legend columns=-1,grid=major,xmode=log,xlabel={$\eta_{\min}$},ylabel={threshold upper endpoint}]
\addplot[blue,mark size=2.9pt,only marks,mark=*] coordinates {(0.000895218714972,1.42626953125) (8.99516713225e-06,1.38000488281) (8.99951616512e-08,1.37548828125) (8.99995220058e-10,1.37475585938)};
\addlegendentry{$p=2$}
\addplot[blue,mark size=2.9pt,only marks,mark=square*,forget plot] coordinates {(0.000155285631915,1.45275878906) (8.99516713253e-06,1.43115234375) (8.99951615729e-08,1.42919921875) (8.99994828663e-10,1.42907714844)};
\addplot[blue,mark size=2.9pt,only marks,mark=triangle*,forget plot] coordinates {(4.18446220432e-05,1.57995605469) (7.49921615693e-06,1.91320800781) (7.49992163998e-08,1.88562011719) (7.49999062855e-10,1.88269042969)};
\addplot[red,mark size=2.2pt,only marks,mark=*] coordinates {(0.000895218714972,1.12573242188) (8.99516713225e-06,0.925720214844) (8.99951616512e-08,0.896240234375) (8.99995220058e-10,0.896240234375)};
\addlegendentry{$p=4$}
\addplot[red,mark size=2.2pt,only marks,mark=*,mark options={fill=white},forget plot] coordinates {(0.000895218714972,0.650756835938) (8.99995220058e-10,0.575134277344)};
\addplot[red,mark size=2.2pt,only marks,mark=square*,forget plot] coordinates {(0.000155285631915,0.957824707031) (8.99516713253e-06,0.956298828125) (8.99951615729e-08,0.956115722656) (8.99994828663e-10,0.956115722656)};
\addplot[red,mark size=2.2pt,only marks,mark=square*,mark options={fill=white},forget plot] coordinates {(0.000155285631915,0.588317871094) (8.99994828663e-10,0.587158203125)};
\addplot[red,mark size=2.2pt,only marks,mark=triangle*,forget plot] coordinates {(4.18446220432e-05,1.10198974609) (7.49921615693e-06,1.37475585938) (7.49992163998e-08,1.3896484375) (7.49999062855e-10,1.39111328125)};
\addplot[red,mark size=2.2pt,only marks,mark=triangle*,mark options={fill=white},forget plot] coordinates {(4.18446220432e-05,0.663635253906) (7.49999062855e-10,0.778564453125)};
\addplot[teal,mark size=1.5pt,only marks,mark=*] coordinates {(0.000895218714972,0.857421875) (8.99516713225e-06,0.997436523438) (8.99951616512e-08,0.790100097656) (8.99995220058e-10,0.777709960938)};
\addlegendentry{$p=6$}
\addplot[teal,mark size=1.5pt,only marks,mark=*,mark options={fill=white},forget plot] coordinates {(0.000895218714972,0.531616210938) (8.99995220058e-10,0.520690917969)};
\addplot[teal,mark size=1.5pt,only marks,mark=square*,forget plot] coordinates {(0.000155285631915,0.800903320312) (8.99516713253e-06,0.802551269531) (8.99951615729e-08,0.802734375) (8.99994828663e-10,0.802734375)};
\addplot[teal,mark size=1.5pt,only marks,mark=square*,mark options={fill=white},forget plot] coordinates {(0.000155285631915,0.486145019531) (8.99994828663e-10,0.484008789062)};
\addplot[teal,mark size=1.5pt,only marks,mark=triangle*,forget plot] coordinates {(4.18446220432e-05,0.873901367188) (7.49921615693e-06,0.98876953125) (7.49992163998e-08,0.990234375) (7.49999062855e-10,0.990478515625)};
\addplot[teal,mark size=1.5pt,only marks,mark=triangle*,mark options={fill=white},forget plot] coordinates {(4.18446220432e-05,0.526672363281) (7.49999062855e-10,0.598876953125)};
\end{axis}\end{tikzpicture}
\caption{Color and decreasing marker size denote $p=2,4,6$; circles, squares and triangles denote axis, diagonal and oblique paths. Filled/open markers use $L=2/4$. Points are not joined across changes in the active mesh or patches.}
\label{fig:sequence-cut}
\end{figure}

With boundary value mass \(M_\partial\) and consistency matrix \(C\), set
\(E=K+0.1S+(p^2/h)M_\partial\) and
\(A(\beta_0)=K+0.1S-C-C^T+\beta_0(p^2/h)M_\partial\).
The nominal coercivity is the smallest energy-normalized stiffness
eigenvalue \(\lambda_{\min}(A(16),E)\). The reported penalty brackets
enclose its crossing of \(10^{-8}\), a small positive numerical margin
above loss of definiteness, rather than the exact zero crossing.
Figure~\ref{fig:sequence-cut} and Table~\ref{tab:sequence-cut} show bounded
thresholds down to \(\eta_{\min}=7.50\times10^{-10}\).
The largest axis-path threshold decreases from 1.4263 to 1.1257 and
0.9974 at \(p=2,4,6\): the \(p^2\) scaling is conservative on this sample.

\begin{table}[!htbp]
\caption{Group minima of lower endpoints and nominal coercivity, and maxima of upper endpoints and flux constants. The cut fraction is the minimum over offsets, not constant along a path. Flux maxima use $\gamma=0.1$ for cuts and $\gamma=0$ for the fitted polygon; dashes mean unmeasured four-layer flux constants.}
\label{tab:sequence-cut}
\centering\footnotesize
\begin{tabular}{lrcrrrrr}
\toprule Geometry & $L$ & $p$ & Min. $\eta_{\min}$ & Lower & Upper & Coercivity & Max. $C_\gamma$ \\ \midrule
axis & 2 & 2,4,6 & $9.00\times10^{-10}$ & 0.7776 & 1.4263 & 0.9121 & 1.5374 \\
axis & 4 & 4,6 & $9.00\times10^{-10}$ & 0.5206 & 0.6508 & 0.9596 & -- \\
diagonal & 2 & 2,4,6 & $9.00\times10^{-10}$ & 0.8008 & 1.4528 & 0.9131 & 1.4823 \\
diagonal & 4 & 4,6 & $9.00\times10^{-10}$ & 0.4839 & 0.5883 & 0.9617 & -- \\
oblique & 2 & 2,4,6 & $7.50\times10^{-10}$ & 0.8738 & 1.9132 & 0.8672 & 2.5694 \\
oblique & 4 & 4,6 & $7.50\times10^{-10}$ & 0.5266 & 0.7786 & 0.9511 & -- \\
fitted & -- & 1--8 & 1 & 1.2957 & 1.7468 & 0.8954 & 1.7696 \\
\bottomrule\end{tabular}
\end{table}

The fitted control uses the same local basis and nominal \(h=1/3\)
on a triangulation of an inscribed regular 20-gon, with no stabilization.
It is a boundary-fitted polygonal approximation, not an exact fitted
representation of the circle. This comparison separates some basis effects from cut effects,
while retaining that geometric difference; it is not an accuracy comparison
for the disk problem. At matching degrees $p=2,4,6$,
the axis and diagonal cuts, including fractions near $10^{-9}$, need no
larger penalty than the fitted control. This is not a universal ordering:
the largest oblique $p=2$ upper endpoint is $1.9132$, versus $1.5232$
for the fitted control at $p=2$.
The evidence is bounded, comparable penalties despite severe cuts.

\FloatBarrier
\subsection{Mesh and Degree Convergence}
\label{subsec:numerical-convergence-2d}
For the smooth experiment, \(u=\exp(x_1+x_2)\), \(f=-2u\), and
\(Eu\) is that same formula on the whole active mesh, fixed independently
of \(h,p\). Alongside \(L^2\) and physical gradient errors, we report
\begin{equation}
E_{\rm phys}^2=\|\nabla(u-u_h)\|_\Omega^2+
\frac{p^2}{h}\|u-u_h\|_{\partial\Omega}^2,
\qquad \energy{Eu-u_h}_h^2=E_{\rm phys}^2+\|Eu-u_h\|_{s_h}^2
\label{eq:physical-error-part}
\end{equation}
Here \(s_h\) includes \(\gamma_1\), and the boundary norm weight excludes
\(\beta_0\). Only the stabilization contribution depends on the chosen
extension outside \(\Omega\).

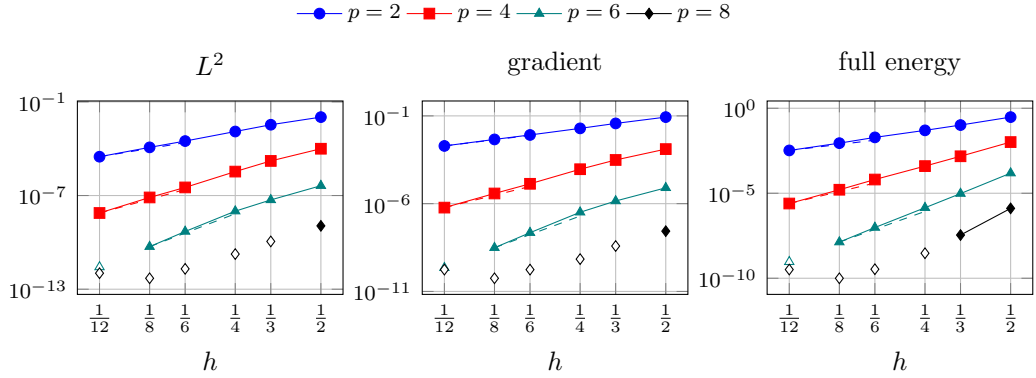
\begin{figure}[!htbp]\centering
\pgfplotslegendfromname{legend-mesh-convergence}\par\vspace{.1em}
\begin{tikzpicture}\begin{axis}[width=.47\textwidth,height=.27\textwidth,font=\small,tick label style={font=\scriptsize},legend style={font=\scriptsize,draw=none},grid=major,width=.32\textwidth,height=.26\textwidth,xmode=log,ymode=log,xlabel={$h$},xtick={0.0833333333333,0.125,0.166666666667,0.25,0.333333333333,0.5},xticklabels={$\frac1{12}$,$\frac18$,$\frac16$,$\frac14$,$\frac13$,$\frac12$},minor x tick num=0,title={$L^2$},legend to name=legend-mesh-convergence,legend columns=-1,legend style={draw=none,fill=none,inner sep=3pt}]
\addplot[blue,mark=*] coordinates {(0.5,0.0105553410103) (0.333333333333,0.00342291532666) (0.25,0.00125887399521) (0.166666666667,0.000307963177218) (0.125,0.000121614151578) (0.0833333333333,3.06840030363e-05)};
\addlegendentry{$p=2$}
\addplot[red,mark=square*] coordinates {(0.5,9.94148271504e-05) (0.333333333333,1.6446457517e-05) (0.25,3.41760503945e-06) (0.166666666667,3.30590597432e-07) (0.125,7.53792056641e-08) (0.0833333333333,7.61089857633e-09)};
\addlegendentry{$p=4$}
\addplot[teal,mark=triangle*] coordinates {(0.5,4.26256049e-07) (0.333333333333,5.23263126814e-08) (0.25,9.97433570054e-09) (0.166666666667,4.97197442905e-10) (0.125,5.29674116034e-11)};
\addlegendentry{$p=6$}
\addplot[teal,mark=triangle*,only marks,mark options={fill=white},forget plot] coordinates {(0.0833333333333,2.63975496134e-12)};
\addplot[black,mark=diamond*] coordinates {(0.5,1.13994462576e-09)};
\addlegendentry{$p=8$}
\addplot[black,mark=diamond*,only marks,mark options={fill=white},forget plot] coordinates {(0.333333333333,1.14956568973e-10) (0.25,1.78220368615e-11) (0.166666666667,2.0118582696e-12) (0.125,5.01538406605e-13) (0.0833333333333,1.05998029666e-12)};
\addplot[blue,dashed,mark=none,forget plot] coordinates {(0.0833333333333,3.06840030363e-05) (0.166666666667,0.00024547202429)};
\addplot[red,dashed,mark=none,forget plot] coordinates {(0.0833333333333,7.61089857633e-09) (0.166666666667,2.43548754443e-07)};
\addplot[teal,dashed,mark=none,forget plot] coordinates {(0.125,5.29674116034e-11) (0.25,6.77982868524e-09)};
\end{axis}\end{tikzpicture}\begin{tikzpicture}\begin{axis}[width=.47\textwidth,height=.27\textwidth,font=\small,tick label style={font=\scriptsize},legend style={font=\scriptsize,draw=none},grid=major,width=.32\textwidth,height=.26\textwidth,xmode=log,ymode=log,xlabel={$h$},xtick={0.0833333333333,0.125,0.166666666667,0.25,0.333333333333,0.5},xticklabels={$\frac1{12}$,$\frac18$,$\frac16$,$\frac14$,$\frac13$,$\frac12$},minor x tick num=0,title={gradient}]
\addplot[blue,mark=*,forget plot] coordinates {(0.5,0.0851609374947) (0.333333333333,0.0371528896088) (0.25,0.0194173219993) (0.166666666667,0.00818222040044) (0.125,0.00455660221145) (0.0833333333333,0.0019392625298)};
\addplot[red,mark=square*,forget plot] coordinates {(0.5,0.00128225292135) (0.333333333333,0.000308104329878) (0.25,9.15034514455e-05) (0.166666666667,1.35274294616e-05) (0.125,3.80044568943e-06) (0.0833333333333,5.95406038422e-07)};
\addplot[teal,mark=triangle*,forget plot] coordinates {(0.5,7.97573398589e-06) (0.333333333333,1.41777287699e-06) (0.25,3.23743356554e-07) (0.166666666667,2.21572920636e-08) (0.125,3.07932262201e-09)};
\addplot[teal,mark=triangle*,only marks,mark options={fill=white},forget plot] coordinates {(0.0833333333333,2.21329308224e-10)};
\addplot[black,mark=diamond*,forget plot] coordinates {(0.5,2.71630691686e-08)};
\addplot[black,mark=diamond*,only marks,mark options={fill=white},forget plot] coordinates {(0.333333333333,3.86131383466e-09) (0.25,6.93666534967e-10) (0.166666666667,1.74567586943e-10) (0.125,5.61722996235e-11) (0.0833333333333,1.74584092806e-10)};
\addplot[blue,dashed,mark=none,forget plot] coordinates {(0.0833333333333,0.0019392625298) (0.166666666667,0.0077570501192)};
\addplot[red,dashed,mark=none,forget plot] coordinates {(0.0833333333333,5.95406038422e-07) (0.166666666667,9.52649661475e-06)};
\addplot[teal,dashed,mark=none,forget plot] coordinates {(0.125,3.07932262201e-09) (0.25,1.97076647809e-07)};
\end{axis}\end{tikzpicture}\begin{tikzpicture}\begin{axis}[width=.47\textwidth,height=.27\textwidth,font=\small,tick label style={font=\scriptsize},legend style={font=\scriptsize,draw=none},grid=major,width=.32\textwidth,height=.26\textwidth,xmode=log,ymode=log,xlabel={$h$},xtick={0.0833333333333,0.125,0.166666666667,0.25,0.333333333333,0.5},xticklabels={$\frac1{12}$,$\frac18$,$\frac16$,$\frac14$,$\frac13$,$\frac12$},minor x tick num=0,title={full energy}]
\addplot[blue,mark=*,forget plot] coordinates {(0.5,0.301483472352) (0.333333333333,0.10258055941) (0.25,0.0500368151855) (0.166666666667,0.0193198167449) (0.125,0.00884206507563) (0.0833333333333,0.00330008532981)};
\addplot[red,mark=square*,forget plot] coordinates {(0.5,0.0102544269729) (0.333333333333,0.00149667169153) (0.25,0.000395346395939) (0.166666666667,6.42004281847e-05) (0.125,1.60828168037e-05) (0.0833333333333,2.48819894586e-06)};
\addplot[teal,mark=triangle*,forget plot] coordinates {(0.5,0.000152556869905) (0.333333333333,9.56675298225e-06) (0.25,1.40462143621e-06) (0.166666666667,9.61196303221e-08) (0.125,1.33714998261e-08)};
\addplot[teal,mark=triangle*,only marks,mark options={fill=white},forget plot] coordinates {(0.0833333333333,9.23394600411e-10)};
\addplot[black,mark=diamond*,forget plot] coordinates {(0.5,1.2886026062e-06) (0.333333333333,3.48115765527e-08)};
\addplot[black,mark=diamond*,only marks,mark options={fill=white},forget plot] coordinates {(0.25,2.96179396653e-09) (0.166666666667,3.4019329177e-10) (0.125,9.96241478678e-11) (0.0833333333333,3.24744111799e-10)};
\addplot[blue,dashed,mark=none,forget plot] coordinates {(0.0833333333333,0.00330008532981) (0.166666666667,0.0132003413192)};
\addplot[red,dashed,mark=none,forget plot] coordinates {(0.0833333333333,2.48819894586e-06) (0.166666666667,3.98111831338e-05)};
\addplot[teal,dashed,mark=none,forget plot] coordinates {(0.125,1.33714998261e-08) (0.25,8.5577598887e-07)};
\end{axis}\end{tikzpicture}
\caption{Smooth mesh refinement on disk A through $N=24$ and $p=8$. Circles, squares, triangles and diamonds denote $p=2,4,6,8$. Open markers fail the norm-specific accuracy checks and are excluded from fits. Dashed guides for $p=2,4,6$ have powers $p+1$ in $L^2$ and $p$ in gradient and energy, normalized at the finest resolved point; they are reference powers, not error bounds.}
\label{fig:sequence-h}
\end{figure}

Figure~\ref{fig:sequence-h} shows disk A.
Mesh slopes agree with the predicted
orders in the resolved range (Table~\ref{tab:sequence-h}). Slopes somewhat
above \(p+1\) and \(p\) are compatible with pre-asymptotic behavior for
the exponential solution and are not stronger asymptotic claims.

\begin{table}[!htbp]
\caption{Least-squares slopes on the finest three consecutive resolved levels in each norm. Values above the reference orders are pre-asymptotic observations, not stronger asymptotic rates.}
\label{tab:sequence-h}
\centering\small
\begin{tabular}{crlrrr}
\toprule Disk & $p$ & Resolved $N$ window & $L^2$ & Gradient & Energy \\ \midrule
A & 2 & 12,16,24 & 3.332 & 2.079 & 2.542 \\
A & 4 & 12,16,24 & 5.455 & 4.510 & 4.684 \\
A & 6 & 8,12,16 & 7.546 & 6.709 & 6.708 \\
B & 2 & 12,16,24 & 3.435 & 2.124 & 2.512 \\
B & 4 & 12,16,24 & 5.695 & 4.704 & 4.686 \\
B & 6 & 8,12,16 & 7.192 & 6.463 & 6.632 \\
\bottomrule\end{tabular}

\end{table}
At \(N=6\), the two disks have nearly coincident exponential degree
curves. Least-squares fits over the resolved degrees give per-degree
factors of approximately 0.06 in \(L^2\) (\(p=1\)--7) and 0.09
in energy (\(p=1\)--8).
Corollary~\ref{cor:analytic-convergence} provides a sufficient Taylor-based
bound, not a sharp prediction of these factors; the entire exponential
extension admits approximation estimates much sharper than a fixed
Taylor-neighborhood bound.

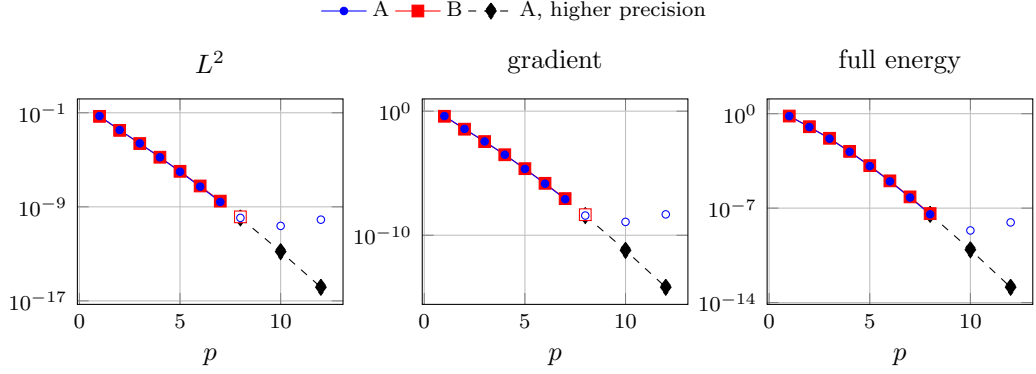
\begin{figure}[!htbp]\centering
\pgfplotslegendfromname{legend-degree-convergence}\par\vspace{.1em}
\begin{tikzpicture}\begin{axis}[width=.47\textwidth,height=.27\textwidth,font=\small,tick label style={font=\scriptsize},legend style={font=\scriptsize,draw=none},grid=major,width=.32\textwidth,height=.27\textwidth,ymode=log,xlabel={$p$},title={$L^2$},legend to name=legend-degree-convergence,legend columns=-1,legend style={draw=none,fill=none,inner sep=3pt}]
\addplot[black,dashed,mark=diamond*,forget plot,mark size=3pt] coordinates {(8,1.15027073493e-10) (10,1.54940830757e-13) (12,1.44226647754e-16)};
\addplot[red,mark=square*,forget plot,mark size=2.2pt] coordinates {(1,0.0496231821795) (2,0.00319308488448) (3,0.000241579783258) (4,1.65909818403e-05) (5,1.05099602104e-06) (6,5.95374713393e-08) (7,3.05747288151e-09)};
\addplot[red,mark=square*,only marks,mark options={fill=white},forget plot,mark size=2.2pt] coordinates {(8,1.4017029904e-10)};
\addplot[blue,mark=*,forget plot,mark size=1.4pt] coordinates {(1,0.0526641757626) (2,0.00342291532666) (3,0.000249564892697) (4,1.6446457517e-05) (5,9.62714894524e-07) (6,5.23263126814e-08) (7,2.61113923214e-09)};
\addplot[blue,mark=*,only marks,mark options={fill=white},forget plot,mark size=1.4pt] coordinates {(8,1.14956568973e-10)};
\addplot[blue,only marks,mark=*,mark options={fill=white},forget plot,mark size=1.4pt] coordinates {(10,2.39291255949e-11) (12,8.10141352645e-11)};
\addlegendimage{blue,mark=*,mark size=1.4pt}
\addlegendentry{A}
\addlegendimage{red,mark=square*,mark size=2.2pt}
\addlegendentry{B}
\addlegendimage{black,dashed,mark=diamond*,mark size=3pt}
\addlegendentry{A, higher precision}
\end{axis}\end{tikzpicture}\begin{tikzpicture}\begin{axis}[width=.47\textwidth,height=.27\textwidth,font=\small,tick label style={font=\scriptsize},legend style={font=\scriptsize,draw=none},grid=major,width=.32\textwidth,height=.27\textwidth,ymode=log,xlabel={$p$},title={gradient}]
\addplot[black,dashed,mark=diamond*,forget plot,mark size=3pt] coordinates {(8,3.87080501513e-09) (10,6.0520056285e-12) (12,6.42104119405e-15)};
\addplot[red,mark=square*,forget plot,mark size=2.2pt] coordinates {(1,0.386090545627) (2,0.0352394652723) (3,0.00351714607448) (4,0.000306005150809) (5,2.31567314108e-05) (6,1.53062594878e-06) (7,8.92314941972e-08)};
\addplot[red,mark=square*,only marks,mark options={fill=white},forget plot,mark size=2.2pt] coordinates {(8,4.52674769652e-09)};
\addplot[blue,mark=*,forget plot,mark size=1.4pt] coordinates {(1,0.4032515945) (2,0.0371528896088) (3,0.00360934709447) (4,0.000308104329878) (5,2.22524337461e-05) (6,1.41777287699e-06) (7,8.0294838728e-08)};
\addplot[blue,mark=*,only marks,mark options={fill=white},forget plot,mark size=1.4pt] coordinates {(8,3.86131383466e-09)};
\addplot[blue,only marks,mark=*,mark options={fill=white},forget plot,mark size=1.4pt] coordinates {(10,1.18716816995e-09) (12,4.80287838993e-09)};
\end{axis}\end{tikzpicture}\begin{tikzpicture}\begin{axis}[width=.47\textwidth,height=.27\textwidth,font=\small,tick label style={font=\scriptsize},legend style={font=\scriptsize,draw=none},grid=major,width=.32\textwidth,height=.27\textwidth,ymode=log,xlabel={$p$},title={full energy}]
\addplot[black,dashed,mark=diamond*,forget plot,mark size=3pt] coordinates {(8,3.48251295334e-08) (10,8.20979069736e-11) (12,1.3591971633e-13)};
\addplot[red,mark=square*,forget plot,mark size=2.2pt] coordinates {(1,0.670368185302) (2,0.108353403558) (3,0.0150578559443) (4,0.00161093655109) (5,0.000141182505133) (6,1.05081456286e-05) (7,6.81153043142e-07) (8,3.9328309299e-08)};
\addplot[blue,mark=*,forget plot,mark size=1.4pt] coordinates {(1,0.659600865927) (2,0.10258055941) (3,0.0141224984351) (4,0.00149667169153) (5,0.000129683941713) (6,9.56675298225e-06) (7,6.09834595811e-07) (8,3.48115765527e-08)};
\addplot[blue,only marks,mark=*,mark options={fill=white},forget plot,mark size=1.4pt] coordinates {(10,2.20569111894e-09) (12,8.72107704717e-09)};
\end{axis}\end{tikzpicture}
\caption{Degree convergence at $N=6$. The higher-precision curve through $p=12$ uses analytic assembly, iterative refinement and coefficient pairs that retain additional precision (Appendix~\ref{sec:numerical-details}). Disks A and B nearly coincide. Open markers retain ordinary-precision failures, including the degree-ten/twelve comparison. The higher-precision points use independently checked coefficient-field norms.}
\label{fig:sequence-p}
\end{figure}

\begin{table}[!htbp]
\caption{Physical and full energy errors for resolved disk-A cases. Shares refer to squared energy.}
\label{tab:physical-energy}
\centering\small
\begin{tabular}{rrrrr}
\toprule $N$ & $p$ & $E_{\rm phys}$ & Full energy & Stabilization share \\ \midrule
6 & 2 & $3.74\times10^{-2}$ & $1.03\times10^{-1}$ & 86.7\% \\
24 & 2 & $1.94\times10^{-3}$ & $3.30\times10^{-3}$ & 65.3\% \\
6 & 4 & $3.09\times10^{-4}$ & $1.50\times10^{-3}$ & 95.7\% \\
24 & 4 & $5.97\times10^{-7}$ & $2.49\times10^{-6}$ & 94.2\% \\
6 & 6 & $1.42\times10^{-6}$ & $9.57\times10^{-6}$ & 97.8\% \\
16 & 6 & $3.08\times10^{-9}$ & $1.34\times10^{-8}$ & 94.7\% \\
\bottomrule\end{tabular}

\end{table}

\begin{figure}[!htbp]\centering
\begin{tikzpicture}\begin{axis}[width=.47\textwidth,height=.27\textwidth,font=\small,tick label style={font=\scriptsize},legend style={font=\scriptsize,draw=none,fill=none,at={(0.5,1.06)},anchor=south,/tikz/column sep=3pt},legend columns=-1,grid=major,ymode=log,xlabel={$p$},ylabel={stabilization seminorm}]
\addplot[red,dashed,mark=square*,forget plot,mark size=2.2pt] coordinates {(1,0.518978834295) (2,0.0955329189988) (3,0.0136506376583) (4,0.00146447952419) (5,0.000127754224467) (6,9.4608356838e-06) (7,6.045139983e-07) (8,3.45963720259e-08)};
\addplot[blue,mark=*,forget plot,mark size=1.4pt] coordinates {(1,0.648071319917) (2,0.11878036414) (3,0.0160686100149) (4,0.00165952571445) (5,0.000140598326132) (6,1.01740185986e-05) (7,6.37908341236e-07) (8,3.58925741102e-08)};
\addlegendimage{blue,mark=*,mark size=1.4pt}
\addlegendentry{$\|Eu\|_{s_h}$}
\addlegendimage{red,dashed,mark=square*,mark size=2.2pt}
\addlegendentry{$\|Eu-u_h\|_{s_h}$}
\end{axis}\end{tikzpicture}
\caption{Extension approximation and error stabilization at $N=6$ on disk A, with the same fixed $Eu$ and $\gamma_1=0.1$.}
\label{fig:extension-stabilization}
\end{figure}
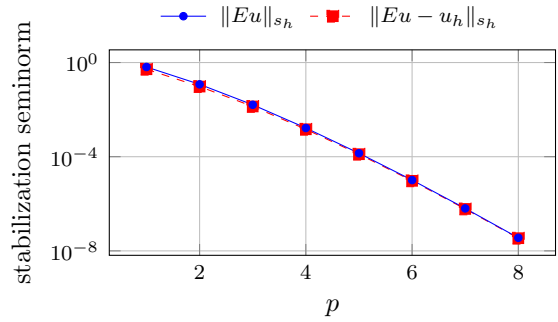
Table~\ref{tab:physical-energy} separates physical and full energy errors, while Figure~\ref{fig:extension-stabilization} compares the extension and error seminorms.
The large stabilization share is consistent with approximation of
\(Eu\) by one polynomial over an entire patch, whose diameter is several
element widths, rather than elementwise approximation alone. The share does not itself
measure physical-accuracy loss; no same-domain fitted or face-jump comparison
establishes comparable accuracy constants here.
At \(p=8\), \(\|Eu\|_{s_h}\) is only 3.7\% (A) or 4.2\% (B)
larger than \(\|Eu-u_h\|_{s_h}\). This close agreement identifies that
extension-approximation scale; it does not turn full energy into a
physical error independent of the extension.

\begin{figure}[!htbp]\centering
\pgfplotslegendfromname{legend-radial-convergence}\par\vspace{.1em}
\begin{tikzpicture}\begin{axis}[font=\small,tick label style={font=\scriptsize},grid=major,width=.32\textwidth,height=.27\textwidth,xmode=log,ymode=log,xtick={2,4,8,16},xticklabels={2,4,8,16},xlabel={$p$},title={$L^2$},legend to name=legend-radial-convergence,legend columns=-1,legend style={font=\scriptsize,draw=none,fill=none,inner sep=3pt}]
\addplot[black,dashed,forget plot] coordinates {(4,0.000940230212567) (8,8.31053948978e-05) (16,7.34554853568e-06)};
\addplot[black,dashdotted,line width=0.8pt,forget plot] coordinates {(4,0.000470115106284) (8,8.31053948978e-05) (16,1.46910970714e-05)};
\addplot[red,mark=square*,forget plot,mark size=2.2pt] coordinates {(2,0.00592927605564) (3,0.00369012808827) (4,0.000722653584863) (5,0.000509131647081) (6,0.000198796831244) (7,0.000155747052902) (8,8.00249093244e-05) (10,3.825198979e-05) (12,2.07672149788e-05)};
\addplot[red,only marks,mark=square*,mark options={fill=white},forget plot,mark size=2.2pt] coordinates {(14,1.24337399745e-05)};
\addplot[blue,mark=*,forget plot,mark size=1.4pt] coordinates {(2,0.00603894255218) (3,0.00378116409459) (4,0.000738517743301) (5,0.00052578479398) (6,0.000207160262825) (7,0.000161046553449) (8,8.31053948978e-05) (10,4.00992272049e-05) (12,2.18938730543e-05)};
\addplot[blue,only marks,mark=*,mark options={fill=white},forget plot,mark size=1.4pt] coordinates {(14,1.30170081925e-05)};
\addlegendimage{blue,mark=*,mark size=1.4pt}
\addlegendentry{A}
\addlegendimage{red,mark=square*,mark size=2.2pt}
\addlegendentry{B}
\addlegendimage{black,dashed}
\addlegendentry{approximation guide}
\addlegendimage{black,dashdotted,line width=0.8pt}
\addlegendentry{theorem guide}
\end{axis}\end{tikzpicture}\begin{tikzpicture}\begin{axis}[font=\small,tick label style={font=\scriptsize},grid=major,width=.32\textwidth,height=.27\textwidth,xmode=log,ymode=log,xtick={2,4,8,16},xticklabels={2,4,8,16},xlabel={$p$},title={physical gradient}]
\addplot[black,dashed,forget plot] coordinates {(4,0.00996531691045) (8,0.00176163579101) (16,0.000311416153452)};
\addplot[black,dashdotted,line width=0.8pt,forget plot] coordinates {(4,0.00704654316405) (8,0.00176163579101) (16,0.000440408947753)};
\addplot[red,mark=square*,forget plot,mark size=2.2pt] coordinates {(2,0.0584474517933) (3,0.0337850232403) (4,0.00970116907161) (5,0.00686759568996) (6,0.00349112131013) (7,0.00274057770069) (8,0.00172950586694) (10,0.000984956206254) (12,0.000622347283185)};
\addplot[red,only marks,mark=square*,mark options={fill=white},forget plot,mark size=2.2pt] coordinates {(14,0.000425449134887)};
\addplot[blue,mark=*,forget plot,mark size=1.4pt] coordinates {(2,0.0588689119146) (3,0.0343231857432) (4,0.00976830077637) (5,0.00697625884187) (6,0.00356972480952) (7,0.00279622796868) (8,0.00176163579101) (10,0.00101136401067) (12,0.000640969581404)};
\addplot[blue,only marks,mark=*,mark options={fill=white},forget plot,mark size=1.4pt] coordinates {(14,0.000437636331939)};
\end{axis}\end{tikzpicture}\begin{tikzpicture}\begin{axis}[font=\small,tick label style={font=\scriptsize},grid=major,width=.32\textwidth,height=.27\textwidth,xmode=log,ymode=log,xtick={2,4,8,16},xticklabels={2,4,8,16},xlabel={$p$},title={full energy}]
\addplot[black,dashed,forget plot] coordinates {(4,0.0319689472523) (8,0.00565136484738) (16,0.000999029601636)};
\addplot[black,dashdotted,line width=0.8pt,forget plot] coordinates {(4,0.0226054593895) (8,0.00565136484738) (16,0.00141284121185)};
\addplot[red,mark=square*,forget plot,mark size=2.2pt] coordinates {(2,0.186858377881) (3,0.070802714807) (4,0.0327875997826) (5,0.0143134374085) (6,0.0113296236324) (7,0.00598879047145) (8,0.00529432175405) (10,0.00281064781302) (12,0.00167274022901)};
\addplot[red,only marks,mark=square*,mark options={fill=white},forget plot,mark size=2.2pt] coordinates {(14,0.00113332713922)};
\addplot[blue,mark=*,forget plot,mark size=1.4pt] coordinates {(2,0.185361338146) (3,0.0756554613128) (4,0.0348011317313) (5,0.0158185124771) (6,0.0123685890948) (7,0.00700520696316) (8,0.00565136484738) (10,0.00303782513692) (12,0.00187526307084)};
\addplot[blue,only marks,mark=*,mark options={fill=white},forget plot,mark size=1.4pt] coordinates {(14,0.00130375619571)};
\end{axis}\end{tikzpicture}
\caption{Finite-regularity degree convergence at $N=4$. Dashed/dash-dotted guides are $p^{-3.5}/p^{-2.5}$ in $L^2$ and $p^{-2.5}/p^{-2}$ in gradient and full energy, normalized at disk A, $p=8$. Open $p=14$ points fail the field criteria; $p=16$ fails before assembly and has no plotted value. Fits use the resolved sampled degrees $4,5,6,7,8,10,12$.}
\label{fig:sequence-radial}
\end{figure}
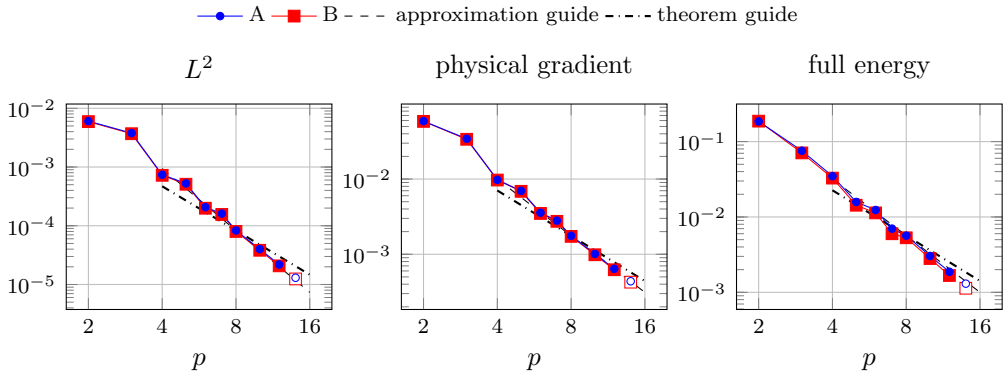

For finite regularity we use \(u=|x-c|^{5/2}\) at \(N=4\),
with the same formula for \(Eu\); Figure~\ref{fig:sequence-radial} shows the degree convergence. Then \(u\in H^s\) for \(s<7/2\),
and the load \(f=-(25/4)|x-c|^{1/2}\) is bounded but nonsmooth
at the center. Both centers lie strictly inside an element; no rate enhancement
from a singularity aligned with a mesh vertex is assumed here.
For each fixed \(0<\varepsilon<1/2\), take \(s=7/2-\varepsilon\) in
Theorems~\ref{thm:hp-energy} and~\ref{thm:l2-error}. On this fixed mesh, the
energy and \(L^2\) bounds are \(O(p^{-2+\varepsilon})\) and
\(O(p^{-5/2+\varepsilon})\), respectively, with constants that may depend on
\(\varepsilon\). The limiting comparison exponents are \(-3.5\)
(optimal \(L^2\) approximation), \(-2.5\) (optimal gradient), and
\(-2\) (energy estimate). They are guides, not assertions of endpoint
Sobolev regularity or sharpness.
Additional paired-rule checks extend the resolved curves to \(p=12\) on both disks.
Over the sampled degrees from $p=4$ through the resolved endpoints, $L^2$ slopes are $-3.33,-3.36$ on A,B, gradient slopes are $-2.56,-2.58$, and energy slopes are $-2.59,-2.61$. The gradient and $L^2$ slopes are within $0.08$ and $0.17$, respectively,
of the limiting approximation orders. Neither the half-power loss of
Theorem~\ref{thm:hp-energy} nor the full-power loss of
Theorem~\ref{thm:l2-error} is visible on this resolved window. These theorems
give upper bounds; this finite range does not establish asymptotic sharpness.

At $p=14$, gradient fields differ between assembly rules by 1.5--6.1\% of their errors, exceeding the 1\% criterion. These points are shown open and excluded from fits; finite-regularity data are not automatically immune to conditioning.
At $p=16$, binary64 bubble-basis construction fails before assembly, whereas independent 80- and 100-digit Gram factorizations succeed. These are arithmetic limits, not counterexamples to the error theorem.

\subsection{Practical Accuracy and Conditioning}
\label{subsec:high-degree-precision}
At errors around \(10^{-10}\) and below, ordinary assembly can fail the
1\% field-discrepancy criterion, consistent with arithmetic limitations;
these points are shown open and excluded from fits.
Accurate assembly and coefficient pairs retaining additional precision at \(N=6\)
recover exponential decay through \(p=12\), with physical \(L^2\) and
gradient errors \(1.44\times10^{-16}\) and \(6.42\times10^{-15}\) on disk A
(Figure~\ref{fig:sequence-p}). This supports an implementation-precision
explanation of the degree-eight stall on that geometry; it does not
validate every flagged finer mesh. The norm measurements use independent
coefficient-field identities; Appendix~\ref{sec:numerical-details} describes
the arithmetic and separate online-evaluation checks.

Table~\ref{tab:sequence-schwarz} separates coefficient conditioning from
energy stability. Four-layer additive Schwarz improves the sampled cut
operators over Jacobi, but the exact interior block precludes a scalable
solver claim. The fitted control is informative on this sample: its Jacobi condition
number rises only from 16.9 at \(p=2\) to 132.6 at \(p=6\), while the
axis-cut value rises from 201 to \(4.65\times10^6\). Thus the hierarchical
bubble basis alone does not explain the severe cut-operator growth;
the geometry and the stabilized cut formulation make a substantial
additional contribution. A different fitted geometry is still a limitation
of this comparison. The axis-cut/fitted ratios are about $12$, $470$ and
$3.5\times10^4$ at $p=2,4,6$. A possible contributor is the
full-element control discussed in Section~\ref{sec:continuation-barrier}:
coefficient norms see the whole active mesh, whereas coercivity needs only
one-sided flux control. This is an interpretation, not a consequence for
$\kappa_J$ of Proposition~\ref{prop:interval-continuation}; basis choice,
diagonal scaling and stabilization also enter, and three degrees do not
establish an exponential asymptotic law.
For the same smallest axis cut at $p=6$, increasing $\gamma_1$ from $0.1$
to $1$ lowers $\kappa_J$ from $4.65\times10^6$ to $9.29\times10^5$;
using $L=4$ at $\gamma_1=0.1$ gives $5.05\times10^5$. Both changes
together give $3.15\times10^5$, still far above the fitted $132.6$.

\begin{table}[!htbp]
\caption{Jacobi-scaled and four-layer Schwarz condition estimates. Cut rows use the smallest offset and stabilization layer count $L=2$; subscript 4 denotes $L_{\rm pc}=4$. The fitted polygon has no stabilization.}
\label{tab:sequence-schwarz}
\centering\small
\begin{tabular}{rrrrrr}
\toprule $p$ & Fitted $\kappa_J$ & Axis $\kappa_J$ & Diag. $\kappa_J$ & Axis $\kappa_{4}$ & Diag. $\kappa_{4}$ \\ \midrule
2 & $1.69\times10^{1}$ & $2.01\times10^{2}$ & $2.78\times10^{2}$ & $3.68\times10^{0}$ & $3.20\times10^{0}$ \\
4 & $6.08\times10^{1}$ & $2.86\times10^{4}$ & $4.46\times10^{4}$ & $2.38\times10^{1}$ & $5.34\times10^{1}$ \\
6 & $1.33\times10^{2}$ & $4.65\times10^{6}$ & $5.07\times10^{6}$ & $3.57\times10^{3}$ & $3.99\times10^{4}$ \\
\bottomrule\end{tabular}

\end{table}

\section{Conclusions}
\label{sec:conclusions}

The exponential cost of full-element polynomial continuation need not enter
symmetric Nitsche stability. A patch projection separates the polynomial,
controlled on inward physical segments, from a residual controlled by the
gradient stabilization on complete elements. Under the stated geometric
mesh restriction, the resulting flux estimate gives the fitted-mesh penalty
scale \(p^2/h\). It also yields a Strang estimate, \(hp\)-energy and
\(L^2\)-error estimates, and exponential degree convergence for analytic
extensions on sufficiently fine meshes.

Translated-disk experiments test normalized flux constants, penalty thresholds
and convergence separately from coefficient conditioning. Unit-stabilization
flux constants lie in \([0.560,0.892]\) for \(p\le8\); unstabilized witnesses
grow consistently with \(\eta_{\min}^{-1/2}\), and matched thresholds remain
below the measured computational flux constants. Accurate assembly resolves
degree decay through \(p=12\). Physical and extension-dependent stabilization
errors are reported separately.

The explicit analytic constant is not certified on these meshes. The fitted
polygon compares stability and conditioning, not disk accuracy. Same-domain
accuracy comparisons and smooth nonconvex tests are outside the scope of
the present study; the dense computations make no scalable-solver claim.

\appendix
\section{Verification Protocol and Patch Geometry}
\label{sec:numerical-details}\label{subsec:two-dimensional-setting}
The numerical companion~\cite{HpNumericalCompanion2026} contains the experiment
code, exact case parameters, reference results with all retained accuracy flags,
and regenerable exact tube certificates. Table~\ref{tab:companion-guide} maps
the studies to the relevant sections of its \emph{Numerical Verification Guide}.
The guide and repository documentation give the corresponding data files,
table identifiers, and commands for reproduction and export.

\begin{table}[!htbp]
\caption{Supporting material in the \emph{Numerical Verification Guide}.
Figure and table numbers refer to this paper; the right column gives guide
section titles. Smooth-case records cover both disks, including the disk-B
values omitted from the main plots.}
\label{tab:companion-guide}
\centering\small
\begin{tabular}{@{}>{\raggedright\arraybackslash}p{.38\textwidth}>{\raggedright\arraybackslash}p{.58\textwidth}@{}}
\toprule
Study in This Paper & Guide Sections \\
\midrule
Finite element basis &
Hierarchical Triangular Basis \\
\addlinespace
Continuation (Figure~\ref{fig:continuation-constants}) &
Continuation Barrier \\
\addlinespace
Flux and matched thresholds (Figures~\ref{fig:direct-flux}--\ref{fig:sequence-cut}, Table~\ref{tab:sequence-cut}) &
Direct Flux Spectra and Matched Thresholds; Cut Stability \\
\addlinespace
Fitted-polygon controls (Tables~\ref{tab:sequence-cut} and~\ref{tab:sequence-schwarz}) &
Fitted Polygon and Extension Seminorms \\
\addlinespace
Smooth convergence (Figures~\ref{fig:sequence-h}--\ref{fig:sequence-p}, Table~\ref{tab:sequence-h}) &
Smooth Mesh and Degree Convergence; Extended Mesh Measurements and Accuracy Flags \\
\addlinespace
Accurate assembly (Figure~\ref{fig:sequence-p}) &
Higher-Degree Accuracy and Precision Limits \\
\addlinespace
Energy components and extension seminorms (Table~\ref{tab:physical-energy}, Figure~\ref{fig:extension-stabilization}) &
Extended Mesh Measurements and Accuracy Flags (paragraph ``Physical Energy Components''); Fitted Polygon and Extension Seminorms \\
\addlinespace
Finite regularity (Figure~\ref{fig:sequence-radial}) &
Finite Regularity in the Full Energy Norm; Radial Degree Extensions and Arithmetic Limits \\
\addlinespace
Conditioning (Table~\ref{tab:sequence-schwarz}) &
A Compact Schwarz Comparison; Schwarz Condition Estimates and Iteration Counts; Stabilization and Coefficient Conditioning \\
\addlinespace
Tube certification (Table~\ref{tab:assignment-audit}) &
Exact Tube Certificates \\
\bottomrule
\end{tabular}
\end{table}

The baseline uses binary64 assembly, QR patch projections and dense solves
with extended-precision iterative refinement. Small residual-Gram eigenvalues
are clipped at the documented tolerance, with discarded contributions recorded.
Two assembly rules must produce fields differing by at most 1\% of the smaller
discretization error in each norm; independently regenerated error rules must
agree to \(10^{-4}\) relatively. Selected smooth cases have additional
reference-field checks. Agreement between two rules cannot exclude shared bias.
All accuracy flags are retained, and unresolved values are excluded from fits.
Flux eigenvalues require agreement between the rules within \(10^{-4}\)
and eigenpair and metric checks within \(10^{-5}\), without regularization.
The companion specifies the quadrature orders, arithmetic and per-case checks,
including the independent radial-load and gradient integrations.

\paragraph{Accurate Assembly and Retained Fields.}
The accurate points in Figure~\ref{fig:sequence-p} use analytic circle moments,
high-precision patch projections and iterative refinement. Each exported
coefficient is stored as a sum of two binary64 numbers to retain additional
precision. Independent coefficient-field identities measure the norms, using
the same 1\% reference-field and \(10^{-4}\) error-rule criteria. No
known-solution stabilization load or diagonal shift is added. Reference fields
are used only for validation, although their assembly shares the underlying
analytic implementation. Separate pointwise checks of solution evaluation do
not certify online arithmetic at every quadrature point. The companion gives
the precision settings, checks, ordinary-precision flags and timings.

\paragraph{Exact Tube Certificates.}\label{subsec:geometry-audit}
An exact full-arc certificate establishes tube depth \(h/16\) on all 30
saved geometries covering the standard PDE study. The seed assignment
fails on 18 of 24 two-layer geometries; reassigning 97 of 5324 open angular
pieces to existing neighboring patches resolves every failure without
changing the method (Remark~\ref{rem:assignment-invariance}).
Table~\ref{tab:assignment-audit} gives the measured maximum overlap and diameter ratios.
The companion also retains the geometry recipes and earlier counterexamples.
Its certificate verifier independently checks the event polynomials, isolated
roots, exact segment coverage and patch membership.
These finite checks do not establish a uniform construction on arbitrary
meshes or untested cuts.

\begin{table}[!htbp]
\caption{Tube certification and measured patch geometry. Four-layer patches need no reassignment.}
\label{tab:assignment-audit}
\centering\small
\begin{tabular}{rrrrrr}
\toprule $L$ & Geometries & Seed passes & Assigned passes & Max. overlap & Max. $\mathrm{diam}(\mathcal P)/h$ \\ \midrule
2 & 24 & 6 & 24 & 9 & $3\sqrt2\simeq4.243$ \\
4 & 6 & 6 & 6 & 18 & $\sqrt{34}\simeq5.831$ \\
\bottomrule\end{tabular}

\end{table}

\paragraph{Explicit Trace and Penalty Constants.}
On the right-isosceles mesh, $\operatorname{diam}(T)=\sqrt2h$ and the
inradius is $\varrho h$, with $\varrho=(2-\sqrt2)/2$. The construction in
Appendix~\ref{app:complete-element-trace} gives
$|(y_0-P)\cdot n(y_0)|\ge\varrho h$. For a circle of radius $R$,
$|n(y)-n(y_0)|=|y-y_0|/R$ and
$|(y-y_0)\cdot n(y_0)|=|y-y_0|^2/(2R)$, whence
\begin{equation}
|(y-P)\cdot n(y)-(y_0-P)\cdot n(y_0)|\le 3h^2/R
\end{equation}
This direct circle calculation replaces $C_{\mathrm{flat}}$ in
Appendix~\ref{app:complete-element-trace} by the sharper constant $3/R$.
Thus $h\le\varrho R/6$ gives $c_*=\varrho/2$. This sufficient mesh restriction
comes from the crude Taylor bound $3h^2/R$ used here; it is not a restriction
of the finite element method. The transverse sign is then constant on the
cut section, and the injectivity argument of
Appendix~\ref{app:complete-element-trace} applies.
Consequently $C_{\mathrm{tr}}=4/\varrho=4(2+\sqrt2)\simeq13.66$. For unit stabilization,
the certified disk depth $c_\delta=1/16$ and conservative Jacobian bound
$c_J=1/2$ give $C_{\mathrm{line}}=32$ and
$\beta_*=4C_{\mathrm{tr}}+256.5=16/\varrho+256.5\simeq311.1$.
For $R=0.62$, this sufficient condition is $h\le\varrho R/6\simeq0.03027$, whereas the
smallest sampled mesh spacing is $1/12$. The calculation therefore
illustrates the analytic chain only on sufficiently fine meshes and
does not certify a numerical trace constant for the present experiments.
The gap to measured thresholds includes the polynomial/residual splitting,
the small certified depth \(h/16\), and the conservative trace and Jacobian
bounds. For the disk, \(J(y,t)=1-t/R\), so the exact lower Jacobian bound
on that tube is \(1-h/(16R)\); using \(c_J=1/2\) adds further slack.
The unit-stabilization value above is a sufficient analytic bound, not
the value used to choose the computational penalty \(\beta_0=16\).

\section{Explicit Constant in the Complete-Element Trace Estimate}
\label{app:complete-element-trace}
The curved-section inverse trace estimate and its quadratic degree dependence
are established in \cite[Lemma 3.1 and Section 6]{WuXiao2019}. We give the
conical argument here to track the geometry constant and the sufficient
mesh restriction used in the explicit penalty bound. No optimality of
the constant \(2d/\varrho\) is asserted.

\begin{proof}[Proof of Lemma~\ref{lem:complete-element-trace}]
Set \(S=\partial\Omega\cap T\); there is nothing to prove if \(S\) has zero surface
measure. We first record the elementary geometric construction used below. Let
\(y_0\in S\) and let \(H_0\) be the tangent hyperplane to \(\partial\Omega\) at
\(y_0\). Shape regularity gives a ball of radius \(\varrho h\) contained in \(T\),
where \(\varrho>0\) is uniform over the mesh family. If \(z_T\) is its
centre, choose between the endpoints \(z_T\pm\varrho h\normal(y_0)\)
of a normal diameter of the closed ball. Their signed distances to \(H_0\)
differ by \(2\varrho h\), so one endpoint \(P\in T\) has distance at least
\(\varrho h\) from \(H_0\).

The \(C^2\) regularity of the boundary implies, uniformly for \(y\in S\),
\begin{equation}
|\normal(y)-\normal(y_0)|\le C_\kappa h,
\qquad
|(y-y_0)\cdot\normal(y_0)|\le C_\kappa h^2 
\label{eq:cut-section-flatness}
\end{equation}
Put \(a=(y_0-P)\cdot\normal(y_0)\) and \(\sigma=\operatorname{sign}(a)\).
Write \(\operatorname{diam}(T)\le C_T h\), so that \(|y-P|\le C_T h\).
Split the difference as
\((y-y_0)\cdot\normal(y_0)+(y-P)\cdot(\normal(y)-\normal(y_0))\).
Then \eqref{eq:cut-section-flatness} gives
\begin{equation}
|(y-P)\cdot\normal(y)-a|\le C_{\mathrm{flat}}h^2,
\qquad C_{\mathrm{flat}}:=C_\kappa(1+C_T)
\end{equation}
uniformly for \(y\in S\). Reduce the geometry-dependent \(h_0\) so that
\(h_0\le\varrho/(2C_{\mathrm{flat}})\).
As \(|a|\ge\varrho h\), we obtain, with \(c_*=\varrho/2\),
\begin{equation}
\sigma(y-P)\cdot\normal(y)\ge c_*h
\qquad \forall y\in S
\label{eq:cone-transversality}
\end{equation}
Thus the transverse sign is the same on all of \(S\), including distinct
connected components.

Consider the conical map
\begin{equation}
\Psi:S\times(0,1)\longrightarrow T,
\qquad
\Psi(y,t)=P+t(y-P)
\label{eq:cut-section-cone-map}
\end{equation}
Its image is contained in \(T\) by convexity, and its volume Jacobian is
\begin{equation}
t^{d-1}|(y-P)\cdot\normal(y)| 
\label{eq:cut-section-cone-jacobian}
\end{equation}
A ray from \(P\) meets \(S\) at most once. Indeed, if it had two intersections,
the segment between them would lie in \(T\) by convexity. Every boundary
intersection on that segment would be transverse by
\eqref{eq:cone-transversality}. Consecutive crossings of \(\partial\Omega\)
alternate between entering and leaving \(\Omega\), so their outward-normal
components along the ray have opposite signs, contradicting the common sign
\(\sigma\). Thus \(\Psi\) is injective. The area formula and
\eqref{eq:cone-transversality} now give, for every nonnegative measurable \(g\),
\begin{equation}
\int_S\int_0^1 t^{d-1}g(\Psi(y,t))\,dt\,dS(y)
\le
\frac{1}{c_*h}\int_T g(x)\,dx 
\label{eq:cone-area-bound}
\end{equation}

For fixed \(y\in S\), every component of
\(t\mapsto\nabla w_h(\Psi(y,t))\) is a polynomial of degree at most \(p-1\).
The weighted estimate \eqref{eq:weighted-endpoint-inverse}, with \(m=p-1\) and
\(\alpha=d-1\), therefore yields
\begin{equation}
|\nabla w_h(y)|^2
\le
p(p+d-1)
\int_0^1 t^{d-1}|\nabla w_h(\Psi(y,t))|^2\,dt 
\label{eq:cone-polynomial-endpoint}
\end{equation}
Since \(|\partial_n w_h|\le|\nabla w_h|\), integration over \(S\), followed by
\eqref{eq:cone-area-bound}, proves
\begin{equation}
h\|\partial_n w_h\|_S^2
\le
\frac{1}{c_*}p(p+d-1)\|\nabla w_h\|_T^2
\le
C_{\mathrm{tr}}p^2\|\nabla w_h\|_T^2 
\end{equation}
where \(p(p+d-1)\le dp^2\) for \(p\ge1\).
This argument applies without change in \(d=2\) and \(d=3\).
\end{proof}

\section{A Physical Cone Alternative}
\label{app:cones}
This conditional variant permits corners; the main theorems use normal tubes.
\begin{lem}[Physical Cone Alternative]
\label{lem:physical-cone}
Let a patch boundary piece \(\Gamma_{\mathcal P}\) admit a finite partition
into Lipschitz surface pieces \(S_j\), with outward normal defined almost
everywhere. Suppose there are points \(a_j\in\mathcal P\cap\Omega\) such that
\(\Psi_j(y,t)=a_j+t(y-a_j)\), for \(y\in S_j\), \(0<t<1\), lies in
\(\mathcal P\cap\Omega\), and
\begin{equation}
|(y-a_j)\cdot\normal(y)|\ge c_c h
\quad\text{for almost every }y\in S_j
\end{equation}
Assume these maps have total volume multiplicity at most \(M_c\), counting
preimages both within and between the pieces. Then every
\(q\in\Pp_p(\mathcal P)\), \(p\ge1\), satisfies
\begin{equation}
h\|\partial_nq\|_{\Gamma_{\mathcal P}}^2
\le\frac{M_c}{c_c}p(p+d-1)
\|\nabla q\|_{\mathcal P\cap\Omega}^2
\le\frac{dM_c}{c_c}p^2\|\nabla q\|_{\mathcal P\cap\Omega}^2
\label{eq:physical-cone-flux}
\end{equation}
\end{lem}

\begin{proof}
Each component of \(\nabla q(\Psi_j(y,t))\) has degree at most \(p-1\)
in \(t\). Apply \eqref{eq:weighted-endpoint-inverse} with \(m=p-1\) and
\(\alpha=d-1\), and sum the components. The volume Jacobian is
\(t^{d-1}|(y-a_j)\cdot\normal(y)|\). Integrating over the pieces and using
the multiplicity bound gives \eqref{eq:physical-cone-flux}.
\end{proof}

\begin{rem}[Scope at Corners]
\label{rem:corner-scope}
Different faces may use different apices. With a cut-independent
complete-element trace estimate and the same residual control, the flux
proof applies with \(C_{\mathrm{cone}}=dM_c/c_c\). Patchwise multiplicity gives
\begin{equation}
C_{\mathrm{inv,cone}}=
\max\{4C_{\mathrm{cone}}N_{\mathcal P},
2C_{\mathrm{tr}}+4C_{\mathrm{cone}}\}
\end{equation}
Cones from different patches may overlap. A global multiplicity bound
$M_c^{\rm glob}$ instead gives
\begin{equation}
C_{\mathrm{inv,cone}}^{\rm glob}
=2C_{\mathrm{tr}}+4dM_c^{\rm glob}/c_c
\end{equation}
without a separate patch-overlap factor. On fixed polygons or polyhedra,
Lemma~\ref{lem:complete-element-trace} applies face by face. Uniform cone
height, multiplicity and patch geometry remain construction requirements.
The main results retain their \(C^2\) assumptions; this stability argument
does not supply dual \(H^2\) regularity at reentrant corners.
\end{rem}

\paragraph{Acknowledgements.}
EB was supported in part by Engineering and Physical Sciences Research
Council grants Nos.~\mbox{EP/V050400/1},
\mbox{EP/T033126/1}, and \mbox{EP/X042650/1}.
PH was supported in part by the Swedish Research Council Grant No.\ 2022-03908.
MGL and KL were supported in part by the Swedish Research Council grants
Nos.\ 2021-04925 and 2025-05562, and the Swedish Research Programme Essence.

\paragraph{Use of AI Tools.}
The initial draft was developed and written by the authors without AI assistance.
During subsequent development and revision, the authors used OpenAI GPT-6 Astra
to assist with reviewing and refining
mathematical arguments, planning numerical experiments, developing and checking
numerical software, preparing figures, and revising the manuscript and
supplementary documentation. The authors directed this work, reviewed and
revised the AI-assisted material, and take full responsibility for the
mathematical arguments, numerical results, and final content.

\paragraph{Data Availability.}
The numerical code, reference data, and verification guide are available at
\url{https://github.com/KarlLarssonUmU/one-sided-hp-cutfem}
and are archived in the numerical companion~\cite{HpNumericalCompanion2026}.

\bigskip
\noindent
\footnotesize {\bf Authors' addresses:}

\smallskip
\noindent
Erik Burman,  \quad \hfill \addressuclshort\\
{\tt e.burman@ucl.ac.uk}

\smallskip
\noindent
Peter Hansbo, \quad \hfill \addressjushort\\
{\tt peter.hansbo@ju.se}

\smallskip
\noindent
Mats G. Larson,  \quad \hfill \addressumushort\\
{\tt mats.larson@umu.se}

\smallskip
\noindent
Karl Larsson, \quad \hfill \addressumushort\\
{\tt karl.larsson@umu.se}

\end{document}